\documentclass[11pt,leqno,twoside]{article}
\def\cvd{\hfill$\Box$}

\def\v{{\rm v}}

\def\int{\mathbb{Z}}
\def\C{\mathbb{C}}
\def\R{\mathbb{R}}

\def\Z{\mathbb{Z}}

\def\Ue{{\cal U}_{\varepsilon}({\mathfrak g})}

\def\im{\text{Im}}

\def\O{{\cal O}}

\def\K{{k}}

\def\proof{\noindent{\bf Proof. }}
\def\Pf{\proof}

\def\pf{\proof}

\def\diag{{\rm diag}}

\def\rk{{\rm rk}}

\def\a{\alpha}

\def\b{\beta}
\def\d{\delta}
\def\l{\lambda}
\def\g{\gamma}

\def\t{\tau}

\def\<#1{\langle #1\rangle}

\def\wJ{w_{\!_J}}

\def\vuoto{\varnothing}
\def\V{{\cal V}}

\usepackage{times}
\usepackage{graphics}
\usepackage{amssymb}
\usepackage{amscd}
\usepackage{amsmath}
\usepackage{fancyhdr}

\input xy
\xyoption{all}

\title{A classification of spherical conjugacy classes}

\newtheorem{theorem}{Theorem}[section]
\newtheorem{lemma}[theorem]{Lemma}
\newtheorem{corollary}[theorem]{Corollary}
\newtheorem{proposition}[theorem]{Proposition}

\newtheorem{remark}[theorem]{Remark}

\newcounter{tigre}
\def\totable{\refstepcounter{tigre}\arabic{tigre}}

\fancypagestyle{plain}{%
\fancyhf{}

}

\author{Mauro Costantini\\
Dipartimento di Matematica Pura ed Applicata\\
Torre Archimede - via Trieste 63 - 35121 Padova - Italy\\
email: costantini@math.unipd.it }
\date{}
\begin{document}
\baselineskip=20pt
\maketitle
\begin{abstract}
  Let $G$ be a simple algebraic group over an
  algebraically closed field $k$. We classify the spherical conjugacy classes of $G$.
\end{abstract}
\section{Introduction}

\newcounter{equat}
\def\theequat{(\arabic{equat})}
\def\equat{\refstepcounter{equat}$$~}
\def\endequat{\leqno{\boldsymbol{(\arabic{equat})}}~$$}

\newcommand{\elem}[1]{\stackrel{#1}{\longto}}
\newcommand{\map}[1]{\stackrel{#1}{\to}}
\def\imp{\Rightarrow}
\def\Imp{\Longrightarrow}
\def\iff{\Leftrightarrow}
\def\Iff{\Longleftrightarrow}
\def\to{\rightarrow}
\def\longto{\longrightarrow}
\def\injto{\hookrightarrow}
\def\rtordu{\rightsquigarrow}

Let $G$ be a simple algebraic group over an algebraically closed field $\K$. In this paper we complete the classification of the spherical conjugacy classes of $G$ (we recall that a conjugacy class $\O$ in $G$ is
called {\it spherical} if a Borel subgroup of $G$ has a
dense orbit on $\O$). There has been a lot of work related to this field, and we refer the reader to the detailed account in the introduction of  \cite{gio3}.
Due to the results in  \cite{gio3} and \cite{mauro-cattiva}, we are left to deal with non-unipotent conjugacy classes either when the characteristic of $k$ is bad, or when $G$ is of type $A_n$ and the characteristic is 2. 

The second goal of this paper is the characterization of spherical conjugacy classes in terms of the {\it dimension formula}, see Theorem \ref{finale}. We finally deduce further consequences of the classification.

\setcounter{equation}{0}
\section{Preliminaries}\label{bepi}

We denote by $\C$ the complex numbers, by $\R$ the reals, by
$\Z$ the integers. 

Let $G$ be a simple algebraic group of rank $n$ over 
$\K$, where $k$ is an algebraically closed field. We fix a maximal torus $T$ of $G$, and a Borel subgroup $B$
 containing $T$:  $B^-$ is the Borel subgroup opposite to $B$, $U$ (respectively $U^-$) is the unipotent radical of $B$ (respectively of $B^-$). Then  $\Phi$ is the set of roots relative
to $T$, and $B$ determines the set of positive roots  
$\Phi^+$, and the simple roots $\Delta=\{\alpha_1,\ldots,\alpha_n\}$, $s_\a$ is the simple reflection associated to $\alpha\in \Phi$. 
We shall use the numbering and the description of the simple roots in terms of the canonical basis $(e_1,\ldots, e_k)$ of an appropriate $\R^k$ as in \cite{bourbaki}, Planches I-IX. For the exceptional groups, we shall write $\b=(m_1,\ldots,m_n)$ for $\b=m_1\a_1+\cdots+m_n\a_n$. 
We identify the Weyl group $W$ with $N/T$, where $N$ is the normalizer of $T$: $w_0$ is the longest element of $W$. 
The real space $E=\R \Phi$  is a Euclidean space, endowed with the $W$-invariant scalar product
$(\a_i,\alpha_j) = d_ia_{ij}$. Here
$\{d_1,\ldots,d_n\}$  are
relatively prime positive integers such that if $D$ is the
diagonal matrix with entries $d_1,\ldots,d_n$, then $DA$ is
symmetric, $A = (a_{ij})$ the 
Cartan matrix. 

We put
$\Pi=\{1,\ldots,n\}$, $\vartheta$ is the symmetry of $\Pi$ induced by $-w_0$. We denote by $\ell$ the usual length function on $W$, and by $rk(1-w)$ the rank of $1-w$ in the geometric representation of $W$.

We use the notation $x_\a(\xi)$, $h_\a(z)$, for $\a\in \Phi$, $\xi\in \K$, $z\in\K^\ast$ as in \cite{yale}, \cite{Carter1}. For $\a\in \Phi$ we put $X_\a=\{x_\a(\xi)\mid \xi\in \K\}$, the root-subgroup corresponding to $\a$, and $H_\a=\{h_\a(z)\mid z\in \K^\ast\}$. 
  Given an element $w\in W$ we shall denote
a representative of $w$ in $N$ by $\dot{w}$. We choose the $x_\a$'s  so that, for all $\a\in \Phi$, $n_\a:=x_\a(1)x_{-\a}(-1)x_\a(1)$
lies in $N$ and has image the reflection $s_\a$ in $W$. Then 
\begin{equation}
x_\a(\xi)x_{-\a}(-\xi^{-1})x_\a(\xi)=h_{\a}(\xi)n_\a\quad,\quad
n_\a^2=h_{\a}(-1)
\label{relazioni}
\end{equation}
$$
n_\a x_\a(x) n_\a^{-1}=x_{-\a}(-x)\quad,\quad
h_\a(\xi)x_\b(x)h_\a(\xi)^{-1}=x_\b(\xi^{\<{\b,\a}}x)
$$
for every $\xi\in\K^\ast$, $x\in k$, $\a$, $\b\in\Phi$, where ${\<{\b,\a}}=\frac{2(\b,\a)}{(\a,\a)}$ (\cite{springer}, Proposition 11.2.1).
The family $(x_\a)_{\a\in\Phi}$ is called a {\it realization} of $\Phi$ in $G$.

We put $T^w=\{t\in T\mid wtw^{-1}=t\}$, $T_2=\{t\in T\mid t^2=1\}$. In particular $T^w=T_2$ if $w=w_0=-1$. We also put $S^w=\{t\in T\mid wtw^{-1}=t^{-1}\}$.

For algebraic groups we use the notation in \cite{Hu2}, \cite{Carter2}. In particular, 
for $J\subseteq \Pi$, $\Delta_J=\{\a_j\mid j\in J\}$, $\Phi_J$ is the corresponding root system, $W_J$ the Weyl group, $P_J$ the standard parabolic subgroup of $G$, $L_J=T\<{X_\a\mid \a\in \Phi_J}$ the standard Levi subgroup of $P_J$.
For $z\in W$ we put 
$
U_z=U\cap z^{-1}U^-z$.
Then the unipotent radical $R_uP_J$ of $P_J$ is $U_{w_0\wJ}$,
where $\wJ$ is the longest element of $W_J$.
Moreover
$
U\cap L_J=U_{\wJ}$
is a maximal unipotent subgroup of $L_J$ (of dimension $\ell(\wJ)$), and $T_J:=T\cap L_J'$ is a maximal torus of $L_J'$.
For unipotent classes in exceptional groups we use the notation in \cite{Carter2}. We use the description of centralizers of involutions as in \cite{Iwa}.

If $X$ is  $G$-variety and $x\in X$, we denote by $G.x$ the $G$-orbit of $x$ and by  $G_x$ the isotropy subgroup of $x$ in $G$. We say that $X$ is {\it spherical} if a Borel subgroup of $G$ has a
dense orbit on $X$. It is well known (\cite{Bri}, \cite{Vin} in characteristic $0$,
\cite{gross}, \cite{knop} in positive characteristic)
 that $X$ is spherical if and only if the set $\cal V$ of $B$-orbits in $X$ if finite.
If $H$ is a closed subgroup of $G$ and the homogeneous space $G/H$ is spherical, we say that $H$ is a spherical subgroup of $G$.

Let $g$ be an element of $G$, with Jordan decomposition $g=su$, $s$ semisimple, $u$ unipotent. Using a terminology slightly different from the usual, we say that $g$ is {\it mixed} if $s\not\in Z(G)$ and $u\not=1$.
For each conjugacy class $\O$ in $G$, $w=w_\O$ is the unique element of $W$ such that $BwB\cap \O$ is open dense in $\O$.

If $x$ is an element of a group $K$ and $H\leq K$, we shall denote by $C(x)$ the centralizer of $x$ in $K$, and by $C_H(x)$ the centralizer of $x$ in $H$. If $x$, $y\in K$, then $x\sim y$ means that $x$, $y$ are conjugate in $K$.

If $H$ is an algebraic group, we denote by $B(H)$ a Borel subgroup of $H$.

We denote the identity matrix of order $k$ by $I_k$.

In the remainder of the paper we shall denote by $p$ the characteristic of $k$ (hence $p$ may be 0).

\setcounter{equation}{0}
\section{The classification}\label{classificazione}

We recall that
the bad primes for the individual types of simple groups are as follows:

none when $G$ has type $A_n$;

$p=2$ when $G$ has type $B_n$, $C_n$, $D_n$;

$p=2$ or $3$ when $G$ has type $G_2$, $F_4$, $E_6$, $E_7$;

$p=2$, $3$ or $5$ when $G$ has type $E_8$.

\medskip
\noindent
For convenience we assume $G$ simply-connected, so that centralizers of semisimple elements are connected. However the classification of spherical conjugacy classes in $G$ is independent of the isogeny class.

We put $\tilde \Pi=\Pi\cup\{0\}$, $\tilde \Delta=\Delta\cup\{\a_0\}$ where $\a_0=-\b$, $\b$ the highest root of $\Phi^+$. Thus $\tilde \Pi$ labels the vertices of the extended Dynkin diagram of the root system $\Phi$. For $J\subset \tilde\Pi$, let $\Phi_J=\Z \{\a_i\mid i\in J\}\cap \Phi$ and
 $$
L_J =\<{T,X_\a\mid\a\in\Phi_J}\>
 $$
This is called a {\it pseudo-Levi subgroup} of $G$ (in the sense of \cite{Som98}). Then the following holds:

 \begin{proposition}\label{sommers} Let $t$ in $G$ be semisimple. Then $C(t)$ is conjugate to a subgroup $L_J$ for some $J\subset \tilde\Pi$. Suppose that the characteristic of $k$ is good for $G$. Let $J\subset \tilde\Pi$. There is $t\in G$ such that $L_J =C(t)$.
 \end{proposition}
 \pf These are \cite[Proposition 30, 32]{mns}.\cvd

We recall some basic facts which have been proved for zero or good, odd characteristic.

\begin{theorem}\label{Giovanna} Let $p\not=2$, and let $\O$ be a spherical conjugacy class of a connected reductive algebraic group. If $\O\cap BwB$ is non-empty, then $w^2=1$.
\end{theorem}
\pf If $p$ is zero, or good and odd this is \cite[Theorem 2.7]{gio1}. The  same proof holds as long as $p\neq2$ (see also \cite[Theorem 2.1]{gio-mauro}).   \cvd

\begin{remark}\label{odd}{\rm 
Let $M(W)$ denote  the Richardson-Springer monoid, i.e. the monoid generated by the symbols $r_\alpha$, for $\alpha\in\Delta$, subject to the braid relations and the relation $r_\alpha^2=r_\alpha$ for $\alpha\in\Delta$. 
Given a spherical $G$-variety,  there is an $M(W)$-action on the set ${\cal V}$ of its $B$-orbits. Under additional conditions, one can also define an action of $W$ on  ${\cal V}$. 
These actions have been introduced in \cite{RS} and \cite{knop}, respectively, and they have been further analyzed and applied in  \cite{closures}, \cite[\S 4.1]{MS}, \cite{springer-schubert}. 
By \cite[Theorem 4.2, b)]{knop}, a case in which the action of $W$ is defined is when $p\not=2$. This allows to extend the proof of \cite[Theorem 2.7]{gio1} to the case when $p\not=2$. We shall come back to this point after the achievement of the classification of spherical conjugacy classes in characteristic 2.} \cvd
\end{remark}

Let $\O$ be a conjugacy class of $G$ and let $\V$ be the set
of $B$-orbits in $\O$. There is a natural map 
$\phi\colon \V\to W$ associating to 
$v\in\V$ the element $w$ in the Weyl group of $G$ for which
$v\subseteq BwB$ (equivalently, for which $v\cap BwB\not=\emptyset$).

\begin{theorem}\label{invo} Let $p\not=2$, and let $\O$ be a conjugacy class in  a connected reductive algebraic group. If ${\rm Im}(\phi)$ contains only involutions in $W$, then $\O$ is spherical.
\end{theorem}
\pf If $p$ is zero, or good and odd this is \cite[Theorem
    5.7]{gio2}. The same proof holds as long as $p\neq2$, once noticed again that the action of $W$ on  ${\cal V}$ is defined.   \cvd

\def\v{{\rm v}}
\def\V{{\cal V}}

\begin{theorem}\label{dime}(\cite[Theorem 25]{CCC}, \cite[Theorem
    4.4]{gio1}) A class $\O$ in a connected reductive algebraic group
      $G$ over an algebraically closed field of zero or good odd characteristic is spherical if and only if there exists $v$ in $\V$ such that $\ell(\phi(v))+{\rm rk}(1-\phi(v))=\dim\O$. If this is the case, $v$ is the dense $B$-orbit in $\O$ and $\phi(v)=w_\O$ (and $v=\O\cap Bw_\O B$).\cvd
\end{theorem} 

The elements $w_\O$ of the Weyl group are involutions, i.e. $w_\O^2=1$, are the unique maximal elements in their conjugacy class and are of the form $w_\O=w_0\wJ$, for a certain $\vartheta$-invariant  subset $J$ of $\Pi$ such that $w_0(\a)=\wJ(\a)$ for every $\a\in \Delta_J$ (\cite[Lemma 3.5]{gio1}, \cite[Corollary 2.11]{chan-lu-to}, \cite{perkins-rowley1}).

\medskip

Strategy of the proof. Let $G_\C$ be the corresponding group over $\C$. We have shown in \cite{CCC} that for every spherical conjugacy class $\cal C$ of $G_\C$ there exists an involution $w=w(\cal C)$ in $W$ such that $\dim {\cal C}=\ell(w)+\rk(1-w)$, with ${\cal C}\cap BwB\not=\emptyset$ (in fact even ${\cal C}\cap BwB\cap B^-\not=\emptyset$). For each group $G$ we introduce a certain set $\O(G)$ of semisimple or mixed conjugacy classes which are candidates for being spherical. For each $\O\in\O(G)$ there is a certain spherical conjugacy class ${\cal C}$ in $G_\C$ such that
$
\dim \O=\dim{\cal C}
$.
Let $w=w_{\cal C}$. Our aim is to show that $\O\cap BwB\not=\emptyset$. 
\medskip

We recall the following result proved in \cite[Theorem 5]{CCC} over $\C$, but which is valid with the same proof over any algebraically closed field.

\begin{proposition}\label{metodo}
Suppose that $\O$ contains an element $x\in B{w}B$. 
Then $$\dim B.x
\geq \ell(w)+\rk(1-w).$$ In
particular $\dim {\O} \geq \ell(w)+\rk(1-w)$.
If, in addition, $\dim\O\leq\ell(w)+\rk (1-w)$ then $\O$ is spherical, 
$w=w_\O$ and $B.x$ is the dense $B$-orbit in $\O$.\cvd
\end{proposition}

If $g$ is in $Z(G)$, then $g\in T$, $\O_g=\{g\}$, $w_\O=1$. In the remaining of the paper we shall consider only non-central conjugacy classes.

We shall use the following result

\begin{lemma}\label{scambio} Assume the positive roots $\b_1,\ldots,\b_\ell$ are such that $[X_{\pm \b_i},X_{\pm \b_j}]=1$ for every $1\leq i<j\leq \ell$. 
Then, for $g=n_{\b_1}\cdots n_{\b_\ell}x_{\b_1}(1)\cdots x_{\b_\ell}(1)$, $h\in T$ such that $\b_i(h)\not=1$ for $i=1,\ldots,\ell$ we have
$$
ghg^{-1}\in BwB\cap B^-
$$
where $w=s_{\b_1}\cdots s_{\b_\ell}$.
\end{lemma}
\pf  For every $i=1,\ldots,\ell$, by (\ref{relazioni}) we have
$$
n_{\b_i}x_{\b_i}(1)\ h \ (n_{\b_i}x_{\b_i}(1))^{-1}=
n_{\b_i}x_{\b_i}(1)\ h  x_{\b_i}(-1)h^{-1}h n_{\b_i}^{-1}=
n_{\b_i}x_{\b_i}(1-\b_i(h))n_{\b_i}^{-1}n_{\b_i} h n_{\b_i}^{-1}=
$$
$$
=x_{-\b_i}(\b_i(h)-1) h_i\in B^-\cap Bs_{\b_i}B
$$
where $h_i=n_{\b_i}h n_{\b_i}^{-1}\in T$.
Hence
$$
ghg^{-1}=
x_{-\b_1}(\b_1(h)-1) \cdots x_{-\b_\ell}(\b_\ell(h)-1) h' \in BwB\cap B^{-}
$$
for a certain $h'\in T$. \cvd

The hypothesis of the Lemma are satisfied for instance if $\b_1,\ldots,\b_\ell$ are pairwise orthogonal and long, as in \cite[Lemma 4.1]{mauro-mathZ}. In characteristic 2, we have
$[X_{\gamma},X_{\delta}]=1$ for every pair $(\gamma,\delta)$ of orthogonal roots.

Let $\O$ be the conjugacy class of $x\in G$. In general the orbit map
$\pi:G/C(x)\to \O$ is a bijective morphism, which may be not separable (i.e. an isomorphism). Nevertheless, we have the following result

\begin{lemma}\label{generale} Let $\O$ be a $G$-orbit with isotropy subgroup $H$. Then $\O$ is spherical if and only if $G/H$ is spherical.
\end{lemma}
\pf This is \cite[Remark 2.14]{FR}.\cvd

\begin{proposition}\label{sovra} Let $g\in G$ with Jordan decomposition
  $g=su$, $s$ semisimple, $u$ unipotent. If $\O_g$ is spherical then $\O_s$ and $\O_u$ are spherical.
\end{proposition}
\pf By Lemma \ref{generale}, $C(g)=C(s)\cap C(u)$ is a spherical subgroup of $G$. Hence both $C(s)$ and $C(u)$ are spherical subgroups of $G$ and, by Lemma \ref{generale}, $\O_s$ and $\O_u$ are spherical.\phantom{abcde}
\cvd
\medskip

For $J\subseteq \Pi$ we put $T_J=T\cap L'_J$, a maximal torus of the derived subgroup $L'_J$ of the standard Levi subgroup $L_J$, so that $B_J:=T_JU_{\wJ}$ is a Borel subgroup of $L'_J$.

\begin{lemma}\label{estensione} Let $\O$ be a conjugacy class of $G$, ${\cal F}\subseteq \O$. Assume there exist $J\subseteq \Pi$ such that ${\cal F}\subseteq L_J$ and $(B_J.x)_{ x\in \cal F}$ is a family of pairwise distinct $B_J$-orbits. Then the family $(B.x)_{ x\in \cal F}$ consists of pairwise distinct $B$-orbits.
 \end{lemma}
 \pf Let $x$, $y$ be elements of $\cal F$, and assume $B.x=B.y$. Then there exists $b\in B$ such that $bxb^{-1}=y$, i.e. $bx=yb$. Since $B=TU_{\wJ}U_{w_0\wJ}$
where  
 $U_{w_0\wJ}$ is the unipotent radical of the standard parabolic subgroup $P_J$, we can write $b=tu_1u_2$, where $t\in T$, $u_1\in U_{\wJ}$ and $u_2\in U_{w_0\wJ}$, so that $tu_1u_2x=ytu_1u_2$. Therefore we get  $tu_1x=ytu_1$. We may decompose $T=T_JS$, where $S=\left(\cap_{i\in J}\ker \a_i\right)^\circ$, $t=t_1t_2$, with $t_1\in T_J$, $t_2\in S$. Then $S\leq C(L_J)$, so that $t_1u_1x=yt_1u_1$. But $t_1u_1$ lies in $B_J$, and we conclude that $B_J.x=B_J.y$. Therefore $x=y$ and we are done.\cvd
 
 \begin{lemma}\label{good} Let $x$ be a semisimple element of $G$ with
 $C(x)=L_J$, a pseudo-Levi subgroup of $G$, and assume $\O_x$ is spherical. Let $\tilde x$ be a semisimple element in $G_{\C}$ such that $C(\tilde x)=L_J$ (in $G_\C$). Then $\O_{\tilde x}$ is spherical.
  \end{lemma}
 \pf First we note that such an $\tilde x$ exists, by Proposition \ref{sommers}.
 By Lemma \ref{generale} and  \cite[Theorem 2.2 (i)]{Brundan}, it follows that $\O_{\tilde x}$ is a spherical semisimple conjugacy class in $G_{\C}$.\cvd

\subsection{Type $A_n$, $n\geq 1$}

For every $i=1,\ldots,\left[\frac{n+1}{2}\right]$, we denote by $\b_i$ the root $e_{i}-e_{n+2-i}$.

\begin{proposition}\label{SL(2)} Let $G=SL(2)$, any characteristic. Let $\O$ be a conjugacy class of $G$. Then $\O\cap Bw_\O B\cap B^-$ is non-empty, $\dim\O=\ell(w_\O)+rk(1-w_\O)$ and $\O$ is spherical.
\end{proposition}
\pf
We may work (and we shall usually do) up to a central element, hence we may assume $\O=\O_x$, $x$ either unipotent or semisimple. If $x$ is unipotent then the result follows from \cite[Proposition 11]{CCC}, whose proof is characteristic-free.
If $x$ is semisimple, then either $x$ is central, or $x$ is regular. In the first case $C(x)=G$, and in the second case we may assume $C(x)=T$.
$$
x=
\left(
\begin{array}{cc}
 f &0 \\
 0 & 1/f \\
\end{array}
\right)
$$
for a certain $f\not=\pm1$. Let
$$
g=
\left(
\begin{array}{cc}
0&1 \\
 -1 & 0\\
\end{array}
\right)
\left(
\begin{array}{cc}
1&1 \\
0& 1\\
\end{array}
\right)
$$
Then
$$
gxg^{-1}\in \O\cap Bw B\cap B^-
$$
where $w=w_0$, with 
$\dim\O=2=\ell(w_0)+rk(1-w_0)$. We conclude by Proposition \ref{metodo}.\cvd

\begin{lemma}\label{non-sferica} Let $H$ be connected and reductive, any characteristic. Then $H$ has a regular spherical conjugacy class if and only if the semisimple part of $H$ is of type $A_1^r$. In this case every conjugacy class is spherical.
\end{lemma}
\pf
Without loss of generality we may assume $H=Z\times G_1\times\cdots\times G_r$, where $Z=Z(H)^\circ$ and $G_i$ is simple for each $i=1,\ldots,r$. Let $n_i=\rk G_i$, $N_i$ the number of positive roots of $G_i$ for $i=1,\ldots,r$. Let $x=(z,x_1,\ldots,x_r)$ be an element of $H$, ${\cal O=\cal O}_x$. Then $\O$ is spherical if and only if each $G_i.x_i$ is spherical in $G_i$, and $x$ is regular if and only if each $x_i$ is regular in $G_i$. Moreover, a spherical $G_i$-conjugacy class in $G_i$ has
 dimension at most $n_i+N_i$, while $G_i.x_i$ is regular in $G_i$ if and only if its dimension is $2N_i$. 
 
If the semisimple part of $H$ is of type $A_1^r$, then every conjugacy class of $H$ is spherical by  Proposition \ref{SL(2)}.

Suppose there exists a regular spherical conjugacy class. Then $2N_i\leq N_i+n_i$ for every $i$, which is possible if and only if $N_i=n_i=1$ for every $i$. Hence the semisimple part of $H$ is of type $A_1^r$.
\cvd

\begin{lemma}\label{GL(3)} Let $H=GL(3)$, any characteristic, $g$ a regular element of $H$. Then there exists a subset
${\cal F}=\{x_m\mid m\in k^\ast\}$ of  $\O_g$ such that $(B(H).x_m)_{ m\in k^\ast}$ consists of pairwise distinct $B(H)$-orbits.
\end{lemma}
\pf See the Appendix.\cvd

\begin{proposition}\label{sferiche SL(n)}
 Let $s$ be a semisimple element of $SL(n+1)$ with at most 2 eigenvalues, any characteristic, $\O$ its conjugacy class. Then $\O\cap Bw_\O B\cap B^-$ is non-empty, and $\O$ is spherical.
\end{proposition}
\pf We work up to a central element, hence we may assume 
$$
s=
\diag(a I_k,b I_{n+1-k})\in SL(n+1)
$$
with $a\not=b$, $1\leq k\leq \left[\frac{n+1}{2}\right]$.
Let
$$
g=n_{\b_1}\cdots n_{\b_k}x_{\b_1}(1)\cdots x_{\b_k}(1)
$$
Then
$$
gsg^{-1}\in \O\cap Bw B\cap B^-
$$
with 
$w=w_{\b_1}\cdots w_{\b_k}$,
$\dim\O=\ell(w)+rk(1-w)$. We conclude by Proposition \ref{metodo} and Lemma \ref{scambio}.\cvd

\begin{theorem}\label{non-sferiche SL(n)}
 Let $g$ be an element of $SL(n+1)$, any characteristic, $g=su$ its Jordan decomposition, $\O$ its conjugacy class. Then $\O$ is spherical if and only if one of the following holds:
\begin{itemize}
\item[(a)] $u=1$ and $s$ has at most 2 eigenvalues;
\item[(b)]  $u\not=1$, $s\in Z(G)$ and $u$ has Jordan blocks of sizes at most 2.
\end{itemize}
\end{theorem}
\Pf
Assume that $\O$ is spherical. Suppose that neither $(a)$ nor $(b)$ hold. Since, by \cite[Theorem 2.2]{knop} every conjugacy class contained in the closure of $\O$ is spherical, without loss of generality we may assume
$$
g=\diag(R,S)\quad ,\quad R\in GL(3), \ S\in GL(n-2), \ S\ \text{diagonal}
$$
with
$$
R=\left(
\begin{array}{ccc}
a & 1 & 0 \\
 0 & a & 1\\
 0&0  & a 
\end{array}
\right)
\text{or}
\left(
\begin{array}{ccc}
a & 1 & 0 \\
 0 & a & 0\\
 0&0  & b 
\end{array}
\right)
\text{or}
\left(
\begin{array}{ccc}
a & 0 & 0 \\
 0 & b & 0\\
 0&0  & c 
\end{array}
\right)
$$
with $a$, $b$, $c$ pairwise distinct. Hence $R$ is regular in $GL(3)$.
Consider the elements
$$
g_m=
\diag(x_m,S)\
$$
for $m\in k^\ast$, where $x_m$ is as defined in Lemma \ref{GL(3)}. We apply Lemma \ref{estensione} with $J=\{1,2\}$:
we put ${\cal F}=\{g_m\mid m\in k^\ast\}\subset L_J$.
The $g_m$'s are all $G$-conjugate to $g$, and pairwise not $B_J$-conjugate. By Lemma \ref{estensione} the family $(B.g_m)_{m\in k^\ast}$ is an infinite family of (distinct) $B$-orbits, a contradiction. Hence either (a) or (b) holds.

The remaining assertions follow by Proposition \ref{sferiche SL(n)}, and from the classification of unipotent classes in zero or odd characteristic (\cite[Theorem 3.2]{gio3} and in characteristic 2 (\cite[Table 1]{mauro-cattiva}).\cvd

\begin{center}
\vskip-20pt
$$
\begin{array}{|c||c|c|c|c|}
\hline
\O  &J& w_\O & C(g)&\dim \O \\ 
\hline
\hline
\begin{array}{c}
\diag(a I_k,b I_{n+1-k})\\
k=1,\ldots,\left[\frac{n+1}{2}\right]\\
a\not=b\\
\end{array}
&J_k&\quad \displaystyle s_{\b_1}\cdots s_{\b_k}\quad& T_1 A_{k-1}A_{n-k}&2k(n+1-k)\\
\hline
\end{array}
$$
\end{center}
\begin{center} Table \totable:  Spherical semisimple classes in $A_n$.
\end{center}
where $w_\O=w_0\wJ$, $J_k=\{k+1,\ldots,n-k\}$, for $k=1,\ldots,\left[\frac{n+1}{2}\right]-1$,
$J_{\left[\frac{n+1}{2}\right]}=\vuoto$. 

\subsection{Type $C_n$ (and $B_n$), $p=2$, $n\geq 2$.}

We put
$\b_i=2e_{i}$ for $i=1,\ldots,n$
 and $\g_i=e_{2i-1}+e_{2i}$ for $\ell=1,\ldots,[\frac n2]$.

\begin{proposition}\label{sferiche Sp(2n)}
 Let $x$ be an element of $Sp(2n)$, any characteristic, $n\geq 2$, $\O$ its conjugacy class. If one of the following holds:
\begin{itemize}
\item[(a)] $x=a_\l={\rm diag}(\lambda I_n,
  \lambda^{-1} I_n)$\quad for $\l\not=\pm1$;
\item[(b)]  $x=c_\l={\rm
  diag}(\lambda, I_{n-1},\lambda^{-1}, I_{n-1})$\quad for $\l\not=\pm1$.
\end{itemize}
Then $\O\cap Bw_\O B\cap B^-$ is non-empty, $\dim\O=\ell(w_\O)+rk(1-w_\O)$, and $\O$ is spherical.
\end{proposition}
\pf 
Assume $
x=
a_\l
$
with $\l\not=\pm1$.
Let
$$
g=n_{\b_1}\cdots n_{\b_n}x_{\b_1}(1)\cdots x_{\b_n}(1)
$$
Then
$$
gxg^{-1}\in \O\cap Bw B\cap B^-
$$
with 
$w=s_{\b_1}\cdots s_{\b_n}=w_0$,
$\dim\O=\ell(w)+rk(1-w)$. 

Assume
$
x=
c_\l
$
with $\l\not=\pm1$.
Let
$$
g=n_{\b_1} n_{\b_2}x_{\g_1}(1)x_{\b_1}(1)
$$
Then
$$
gxg^{-1}\in \O\cap Bw B\cap B^-
$$
with 
$w=s_{\b_1}s_{\b_2}$,
$\dim\O=\ell(w)+rk(1-w)$. 

We conclude by Proposition \ref{metodo}.\cvd

\begin{proposition}\label{semisemplici Sp(2n)}
Let $G=Sp_{2n}(k)$, $p=2$, $n\geq2$. The spherical semisimple
  classes are represented by
  \begin{itemize}
\item[(a)]   $a_\lambda={\rm diag}(\lambda I_n,
  \lambda^{-1} I_n)$ \quad for $\lambda\neq 1$,
\item[(b)]   $c_\lambda={\rm
  diag}(\lambda, I_{n-1},\lambda^{-1}, I_{n-1})$ \quad for $\lambda\neq 1$. 
    \end{itemize}
  \end{proposition}
\pf Let $x$ be a semisimple element of $G$, and assume $\O=\O_x$ is spherical. Without loss of generality $C(x)=L_J$, a pseudo-Levi subgroup of $G$. There exists a semisimple element $\tilde x$ in $G_{\C}$ such that $C(\tilde x)$ is $L_J$ in $G_\C$. By Lemma \ref{good}, it follows that $\O_{\tilde x}$ is a spherical semisimple conjugacy class in $G_{\C}$, and therefore, from the classification of semisimple spherical conjugacy classes in zero (or odd) characteristic (\cite[Table 1]{CCC}, \cite[Theorem 3.3]{gio3}), it follows that
$L_J$ is of type
$C_\ell C_{n-\ell}$ for $\ell=1,\ldots,\left[\frac{n}{2}\right]$, $T_1 C_{n-1}$ or $T_1\tilde A_{n-1}$. But $Z(C_\ell C_{n-\ell})=1$, so that we are left with
$$
 a_\lambda={\rm diag}(\lambda I_n,
  \lambda^{-1} I_n) \longleftrightarrow    T_1\tilde A_{n-1} 
  $$
  $$
  c_\lambda={\rm
  diag}(\lambda, I_{n-1},\lambda^{-1}, I_{n-1}) \longleftrightarrow   T_1 C_{n-1}
  $$
  for $\lambda\neq 1$.
We conclude by Proposition \ref{sferiche Sp(2n)}\cvd

\begin{center}
\vskip-20pt
$$
\begin{array}{|c||c|c|c|c|}
\hline
\O  &J& w_\O & C(g)&\dim \O \\ 
\hline
\hline 
\begin{array}{c}
c_\lambda={\rm
  diag}(\lambda, I_{n-1},\lambda^{-1}, I_{n-1})
\\
\lambda\neq 1
\end{array} &\displaystyle J_2 & s_{\b_1} s_{\b_2} &T_1C_{n-1}&4n-2\\
\hline
\begin{array}{c}
a_\lambda={\rm diag}(\lambda I_n,
  \lambda^{-1} I_n)
\\
\lambda\neq 1
\end{array}
&  \displaystyle \emptyset& \quad\displaystyle w_0=s_{\b_1}\cdots s_{\b_n}\quad& T_1\tilde A_{n-1} &n^2+n
\\
\hline
\end{array}
$$
\end{center}
\begin{center} Table \totable\label{Cn}:  Spherical semisimple classes in $C_n$,
$n\geq 2$, $p=2$.
\end{center}
where $w_\O=w_0\wJ$,  $J_2=\vuoto$ if $n=2$, $J_2=\{3,\ldots,n\}$ if $n\geq 3$.

\medskip

We now deal with mixed conjugacy classes. 

\begin{lemma}\label{Sp(4)} Let $H=Sp(4)$, any characteristic, $g$ a mixed regular element of $H$. 
Then there exists a subset
${\cal F}=\{x_m\mid m\in k^\ast\}$ of  $\O_g$ such that $(B(H).x_m)_{ m\in k^\ast}$ consists of pairwise distinct $B(H)$-orbits.
\end{lemma}
\pf See the Appendix.\cvd

\begin{proposition}\label{mixed Sp(2n)} Let $\O$ be the conjugacy class of a mixed element $g$ of $Sp(2n)$, $p=2$. Then $\O$ is not spherical.
\end{proposition}
\Pf
Let $g=su$, Jordan decomposition. Assume, for a contradiction, that $\O$ is spherical. Then both $\O_s$ and $\O_u$ are spherical.
By Proposition \ref{semisemplici Sp(2n)}, $H=C(s)$ is of type $T_1 C_{n-1}$ or $T_1\tilde A_{n-1}$. However $\dim T_1\tilde A_{n-1}=n^2$,
therefore $C_{T_1\tilde A_{n-1}}(u)$ is not spherical in $G$. We are left with $H$ of type $T_1 C_{n-1}$, and we may assume $s=c_a=h_{\b_{1}}(a)$ for a certain $a\not=1$. 

Since every conjugacy class contained in the closure of $\O$ is spherical, it is enough to deal with the minimal non-trivial spherical unipotent classes in $T_1 C_{n-1}$. From the classification of spherical unipotent classes in characteristic 2 (\cite[Table1, Table 2]{mauro-cattiva}), we may assume:
$$
g=h_{\b_{1}}(a)x_{\a_{2}}(1)\quad\text{if $n=2$}
$$
$$
g=h_{\b_{n-1}}(a)x_{\a_n}(1)\quad  {\text {or}} \quad g=h_{\b_{1}}(a)x_{\a_{2}}(1)\quad\text{if $n\geq 3$}
$$
since $h_{\b_{n-1}}(a)=\diag(I_{n-2},a,1,I_{n-2,}a^{-1},1)$ is conjugate to $h_{\b_{1}}(a)$.

Suppose $g=h_{\b_{n-1}}(a)x_{\a_n}(1)$, $n\geq 2$.

\noindent
We apply Lemma \ref{estensione} with $J=\{n-1,n\}$. By considering the corresponding embedding of $C_2$ into $C_n$, we may assume that the family ${\cal F}=\{g_m\mid m\in k^\ast\}$ is a subset of $L_J$.
The $g_m$'s are all $G$-conjugate to $g$, and pairwise not $B_J$-conjugate. By Lemma \ref{estensione} the family $(B.g_m)_{m\in k^\ast}$ is an infinite family of (distinct) $B$-orbits, a contradiction. Hence the class of $g=h_{\b_{n-1}}(a)x_{\a_n}(1)$ is not spherical.

Suppose $g=h_{\b_{1}}(a)x_{\a_{2}}(1)$, $n\geq 3$.

\noindent
Then
$$
g =
\left(
\begin{array}{cccc}
A& 0 & 0 & 0  \\
 0 &   I_{n-3}& 0 &0 \\
 0& 0 &   {}^tA^{-1} 
  & 0  \\
 0 & 0 & 0 &   I_{n-3} \\
\end{array}
\right)
 \quad,\quad
A=
\left(
\begin{array}{ccc}
 a & 0 &0  \\
 0 & 1 &1\\
 0& 0 & 1
\end{array}
\right)
$$
Let $(x_m)_{m\in k^\ast}$ be the family introduced in Lemma \ref{GL(3)}, such that $x_m$ is $GL(3)$-conjugate to $A$ for every $m\in k^\ast$ and $(B(GL(3)).x_m)_{ m\in k^\ast}$ consists of pairwise distinct $B(GL(3))$-orbits. 

We put
$$
g_m =
\left(
\begin{array}{cccc}
  x_m& 0 & 0 & 0  \\
 0 & I_{n-3} & 0 &0 \\
 0& 0 & {}^tx_m^{-1}& 0  \\
 0 & 0 & 0 & I_{n-3} \\
\end{array}
\right)
$$
The $g_m$'s are all $Sp(2n)$-conjugate to $g$. By Lemma \ref{estensione} with $J=\{1,2\}$,
the family $(B.g_m)_{m\in k^\ast}$ is an infinite family of (distinct) $B$-orbits, a contradiction.
Hence the class of $g=h_{\b_{1}}(a)x_{\a_2}(1)$ is not spherical.\cvd

\begin{theorem}\label{Sp(2n)}
Let $G=Sp_{2n}(k)$, $p=2$, $n\geq2$. The spherical
  classes are either semisimple or unipotent.  The semisimple classes are represented in Table \ref{Cn}, the unipotent classes are represented in \cite[Table 2]{mauro-cattiva}.\cvd
\end{theorem}
\medskip

\subsection{Type $D_n$, $p=2$, $n\geq 4$.}

Let $r=\left[\frac{n}2\right]$. We put $\b_\ell=e_{2\ell-1}+e_{2\ell}$, $\d_\ell=e_{2\ell-1}-e_{2\ell}$ for $\ell=1,\ldots, r$.
We put $J_1=\{3,\ldots,n\}$, 
$K_r=\{1,3,\ldots,2r-1\}$. Also if $n$ is even we put $K_r'=\{1,3,\ldots,n-3,n\}$.

In this section we deal with groups $G$ of type $D_n$. We recall that we are assuming $G$ simply-connected. Since $p=2$, the covering map $\pi:G\to SO(2n)$ is an isomorphism of abstract groups. For convenience of the reader, given an element $x\in G$ we shall also indicate the corresponding $\pi(x)\in SO(2n)$.

\begin{proposition}\label{sferiche D_n}
 Let $x$ be an element of $G=D_n$, any characteristic, $n\geq 4$, $\O$ its conjugacy class. If one of the following holds:
\begin{itemize}
\item[(a)] $x=c_\l=h_{\b_1}(\l)h_{\d_1}(\l)$\quad for $\l\not=\pm1$;
\item[(b)]  $x=a_\l=h_{\b_1}(\l)\cdots h_{\b_r}(\l)$\quad for $\l\not=\pm1$.
\item[(c)]  $x=a_\l'=h_{\b_1}(\l)\cdots h_{\b_{r-1}}(\l)h_{\a_{n-1}}(\l)$\quad for $\l\not=\pm1$,\ $n$ even.
\end{itemize}
Then $\O\cap Bw_\O B\cap B^-$ is non-empty, $\dim\O=\ell(w_\O)+rk(1-w_\O)$, and $\O$ is spherical.
\end{proposition}
\pf 
Assume $
x=
c_\l
$
with $\l\not=\pm1$.
Let
$$
g=n_{\b_1}n_{\d_1}x_{\b_1}(1) x_{\d_1}(1)
$$
Then
$$
gxg^{-1}\in \O\cap Bw B\cap B^-
$$
with 
$w=s_{\b_1}s_{\d_1}$,
$\dim\O=\ell(w)+rk(1-w)$. 

Assume
$
x=
a_\l
$
with $\l\not=\pm1$.
Let
$$
g=n_{\b_1}\cdots n_{\b_r}x_{\b_1}(1)\cdots x_{\b_r}(1)
$$
Then
$$
gxg^{-1}\in \O\cap Bw B\cap B^-
$$
with 
$w=s_{\b_1}\cdots s_{\b_r}$,
$\dim\O=\ell(w)+rk(1-w)$. 

$c)$ follows from $b)$, by using the graph automorphism of $G$ exchanging $n-1$ and $n$.
We conclude by Proposition \ref{metodo}.\cvd

\begin{proposition}\label{semisemplici D_n}
Let $G=D_n$, $p=2$, $n\geq 4$. The spherical semisimple
  classes are represented by
 \begin{itemize}
\item[(a)] $x=c_\l=h_{\b_1}(\l)h_{\d_1}(\l)$\quad for $\l\not=1$;
\item[(b)]  $x=a_\l=h_{\b_1}(\l)\cdots h_{\b_m}(\l)$\quad for $\l\not=1$.
\item[(c)]  $x=a_\l'=h_{\b_1}(\l)\cdots h_{\b_{m-1}}(\l)h_{\a_{n-1}}(\l)$\quad for $\l\not=1$,\ $n$ even.
\end{itemize}
  \end{proposition}
\pf Let $x$ be a semisimple element of $G$, and assume $\O=\O_x$ is spherical. Without loss of generality $C(x)=L_J$, a pseudo-Levi subgroup of $G$. There exists a semisimple element $\tilde x$ in $G_{\C}$ such that $C(\tilde x)$ is $L_J$ in $G_\C$. By Lemma \ref{good}, it follows that $\O_{\tilde x}$ is a spherical semisimple conjugacy class in $G_{\C}$, and therefore, from the classification of semisimple spherical conjugacy classes in zero (or odd) characteristic (\cite[Table 1]{CCC}, \cite[Theorem 3.4]{gio3}), it follows that
$L_J$ is of type
$D_\ell D_{n-\ell}$, $2\leq \ell\leq r$, $T_1 D_{n-1}$ or $T_1A_{n-1}$ (up to graph automorphism). But $Z(D_\ell D_{n-\ell})=1$, $2\leq \ell\leq r$  (since $p=2$), so that we are left with
$$
   T_1 D_{n-1}\longleftrightarrow   c_\l =h_{\b_1}(\l)h_{\a_1}(\l)\rightarrow
   \diag(\l^2,I_{n-1},\l^{-2},I_{n-1})
  $$
  $$
T_1 A_{n-1}
  \longleftrightarrow
  a_\l=h_{\b_1}(\l)\cdots h_{\b_r}(\l)\rightarrow
  \diag(\l I_{n},\l^{-1}I_{n})
  $$
  $$
(T_1 A_{n-1})'
  \longleftrightarrow
 a_\l'=h_{\b_1}(\l)\cdots h_{\b_{r-1}}(\l)h_{\a_{n-1}}(\l)\rightarrow
\diag(\l I_{n-1},\l^{-1},\l^{-1}I_{n-1},\l)\ n\text{ even}
  $$
  for $\lambda\neq 1$.
We conclude by Proposition \ref{sferiche D_n}.\cvd
\vfill\eject
\begin{center}
\vskip-20pt
$$
\begin{array}{|c||c|c|c|c|}
\hline
\O  &J& w_\O & SO(2n)&\dim \O \\ 
\hline
\hline 
\begin{array}{c}
c_\l=h_{\b_1}(\l)h_{\d_1}(\l)
\\
\l\not=1
\end{array} &\displaystyle J_1 & s_{\b_1}s_{\d_1}&\diag(\l^2,I_{n-1},\l^{-2},I_{n-1})&4(n-1)\\
\hline
\begin{array}{c}
a_\l=h_{\b_1}(\l)\cdots h_{\b_r}(\l)
\\
\lambda\neq 1
\end{array}
&  \displaystyle K_{r}& \quad\displaystyle s_{\b_1}\cdots s_{\b_r}\quad&  \diag(\l I_{n},\l^{-1}I_{n})&n^2-n
\\
\hline
\begin{array}{c}
a'_\l=h_{\b_1}(\l)\cdots h_{\b_{r-1}}(\l)h_{\a_{n-1}}(\l)\\
\lambda\neq 1
\end{array}
&  \displaystyle K_r'& \quad\displaystyle s_{\b_1}\cdots s_{\b_{r-1}}s_{\a_{n-1}}\quad& \diag(\l I_{n-1},\l^{-1},\l^{-1}I_{n-1},\l)&n^2-n\\
\hline
\end{array}
$$
\end{center}
\begin{center} Table \totable\label{Dn even}:  Spherical semisimple classes in $D_n$, $p=2$,
$n\geq 4$, $n=2r$. 
\end{center}

\begin{center}
\vskip-20pt
$$
\begin{array}{|c||c|c|c|c|}
\hline
\O  &J& w_\O & SO(2n)&\dim \O \\ 
\hline
\hline 
\begin{array}{c}
c_\l=h_{\b_1}(\l)h_{\d_1}(\l)
\\
\lambda\neq 1
\end{array} &\displaystyle J_1 & s_{\b_1}s_{\d_1}&\diag(\l^2,I_{n-1},\l^{-2},I_{n-1})&4(n-1)\\
\hline
\begin{array}{c}
a_\l=h_{\b_1}(\l)\cdots h_{\b_r}(\l)
\\
\lambda\neq 1
\end{array}
&  \displaystyle K_{r}& \quad\displaystyle s_{\b_1}\cdots s_{\b_r}\quad&  \diag(\l I_{n},\l^{-1}I_{n})&n^2-n
\\
\hline
\end{array}
$$
\end{center}
\begin{center} Table \totable\label{Dn odd}:  Spherical semisimple classes in $D_n$, $p=2$,
$n\geq 5$, $n=2r+1$.
\end{center}
where $w_\O=w_0\wJ$.

We now deal with mixed conjugacy classes.

\begin{proposition}\label{non-sferiche SO(2n)} Let $\O$ be the conjugacy class of a mixed element $g$ in $D_n$, $p=2$. Then $\O$ is not spherical. 
\end{proposition}
\pf 
We work with $SO(2n)$ via $\pi$. Let $g=su$, Jordan decomposition. Assume that $\O$ is spherical. Then both $\O_s$ and $\O_u$ are spherical, and we may assume, up to conjugation and graph automorphism, 
that 
$$
s=\diag(a I_{n-1},a^{-1},a^{-1}I_{n-1},a)
\quad
\text{or}
\quad
s=\diag( I_{n-3},a, I_{2},I_{n-3},a^{-1}, I_{2})
$$
for a certain $a\not=1$.

\noindent
Assume $s=\diag(a I_{n-1},a^{-1},a^{-1}I_{n-1},a)$ for a certain $a\not=1$.

Since every conjugacy class contained in the closure of $\O$ is spherical, without loss of generality we may assume
$
u=x_{\a_{n-2}}(a^{-1})
$, so that
$$
g =
\left(
\begin{array}{cccc}
 a I_{n-3} & 0 & 0 & 0  \\
 0 & A& 0 &0 \\
 0& 0 & a^{-1} I_{n-3}
  & 0  \\
 0 & 0 & 0 &   {}^tA^{-1} \\
\end{array}
\right)
 \quad,\quad
A=
\left(
\begin{array}{ccc}
 a & 1 &0  \\
 0 & a &0\\
 0& 0 & a^{-1} 
\end{array}
\right)
$$
Let $(x_m)_{m\in k^\ast}$ be the family introduced in Lemma \ref{GL(3)}, such that $x_m$ is $GL(3)$-conjugate to $A$ for every $m\in k^\ast$ and $(B(GL(3)).x_m)_{ m\in k^\ast}$ consists of pairwise distinct $B(GL(3))$-orbits. 

We put
$$
g_m =
\left(
\begin{array}{cccc}
  a I_{n-3} & 0 & 0 & 0  \\
 0 & x_m & 0 &0 \\
 0& 0 & a^{-1} I_{n-3}& 0  \\
 0 & 0 & 0 & {}^tx_m^{-1}\\
\end{array}
\right)
$$
The $g_m$'s are all $SO(2n)$-conjugate to $g$. By Lemma \ref{estensione} with $J=\{n-1,n-2\}$,
the family $(B.g_m)_{m\in k^\ast}$ is an infinite family of (distinct) $B$-orbits, a contradiction.
This settles the cases when $n$ is odd and $C(s)$ is of type $T_1A_{n-1}$, and $n$ even and $C(s)$ of type $(T_1A_{n-1})'$. Applying the graph automorphism exchanging $n$ and $n-1$, this also settles the case when $n$ is even and $C(s)$ is of type $T_1A_{n-1}$.

\noindent
Assume $s=\diag( a, I_{n-1},a^{-1}, I_{n-1})
$
for a certain $a\not=1$.

Since every conjugacy class contained in the closure of $\O$ is spherical, without loss of generality we may assume
$
u=x_{\a_{2}}(1)
$, so that
$$
g =
\left(
\begin{array}{cccc}
A& 0 & 0 & 0  \\
 0 &   I_{n-3}& 0 &0 \\
 0& 0 &   {}^tA^{-1} 
  & 0  \\
 0 & 0 & 0 &   I_{n-3} \\
\end{array}
\right)
 \quad,\quad
A=
\left(
\begin{array}{ccc}
 a & 0 &0  \\
 0 & 1 &1\\
 0& 0 & 1
\end{array}
\right)
$$
Let $(x_m)_{m\in k^\ast}$ be the family introduced in Lemma \ref{GL(3)}, such that $x_m$ is $GL(3)$-conjugate to $A$ for every $m\in k^\ast$ and $(B(GL(3)).x_m)_{ m\in k^\ast}$ consists of pairwise distinct $B(GL(3))$-orbits. 

We put
$$
g_m =
\left(
\begin{array}{cccc}
  x_m& 0 & 0 & 0  \\
 0 & I_{n-3} & 0 &0 \\
 0& 0 & {}^tx_m^{-1}& 0  \\
 0 & 0 & 0 & I_{n-3} \\
\end{array}
\right)
$$
The $g_m$'s are all $SO(2n)$-conjugate to $g$. By Lemma \ref{estensione} with $J=\{1,2\}$,
the family $(B.g_m)_{m\in k^\ast}$ is an infinite family of (distinct) $B$-orbits, a contradiction.
This settles the case $C(s)$ is of type $T_1D_{n-1}$, and we are done.\cvd

\begin{theorem}\label{D_n}
Let $G=D_n$, $p=2$, $n\geq 4$. The spherical
  classes are either semisimple or unipotent.  The semisimple classes are represented in Table \ref{Dn even}, \ref{Dn odd}, the unipotent classes are represented in \cite[Table 3, 4]{mauro-cattiva}.\cvd
\end{theorem}

\subsection{Type $E_6$.}  

We put
$$
\begin{array}{cclll}
\b_1 = (1,2,2,3,2,1),&
\b_2 =(1,0,1,1,1,1)\\
\b_3 =(0,0,1,1,1,0),&
\b_4 = (0,0,0,1,0,0)
\end{array}
$$

\begin{proposition}\label{sferiche E_6}
 Let $x$ be an element of $E_6$, any characteristic, $\O$ its conjugacy class. If one of the following holds:
\begin{itemize}
\item[(a)] $x=h_{\a_1}(-1)h_{\a_4}(-1)h_{\a_6}(-1)$;
\item[(b)]  $x=h(z)=h_{\a_1}(z^4)h_{\a_2}(z^3)h_{\a_3}(z^5)h_{\a_4}(z^6)h_{\a_5}(z^4)h_{\a_6}(z^2)$\quad for $z^3\not=1$.
\end{itemize}
Then $\O\cap Bw_\O B\cap B^-$ is non-empty, $\dim\O=\ell(w_\O)+rk(1-w_\O)$, and $\O$ is spherical.
\end{proposition}
\pf 
\noindent $(a)$\quad
If $p=2$, then $x=1$, and there is nothing to prove. So assume $p\not=2$. In $G$ there are two classes of involutions: one has centralizer of type $A_1 A_5$ and dimension 40, the other has centralizer of type $D_5T_1$ and has dimension 32.
Let 
$$
y=n_{\b_1}\cdots n_{\b_4}\in w_0B
$$
$w=s_{\b_1}\cdots s_{\b_4}=w_0$.
Then $y^2=h_{\b_1}(-1)\cdots h_{\b_4}(-1)=1$, and $\dim\O_y\geq 40$ by Proposition \ref{metodo}. Since $C(x)$ is of type $A_1A_5$, we conclude that $x\sim y$, so that $\O\cap Bw_0 B$ is non-empty, $\dim\O=\ell(w_0)+rk(1-w_0)$ and $\O$ is spherical. It is a general fact that if $t$ is semisimple and $\O_t\cap Bw B\not=\emptyset$, then  $\O_t\cap Bw B\cap B^-\not=\emptyset$, \cite[Lemma 14]{CCC}.

\noindent $(b)$\quad
 In this case $C(x)$ is of type $D_5T_1$ (note that $C(x)=C(h(-1))$ if $p\not=2$).
Let
$$
g=n_{\b_1} n_{\b_2}x_{\b_1}(1)x_{\b_2}(1)
$$
Then
$$
gxg^{-1}\in \O\cap Bs_{\b_1}s_{\b_2}B\cap B^-
$$
with 
$w=s_{\b_1}s_{\b_2}$,
$\dim\O=\ell(w)+rk(1-w)$. 
We conclude by Proposition \ref{metodo}.\cvd
  
  \begin{proposition}\label{semisemplici E_6}
Let $G=E_6$. The spherical semisimple
  classes are represented by
  $$
   \begin{array}{lc}
  \\
 h(z)=h_{\a_1}(z^4)h_{\a_2}(z^3)h_{\a_3}(z^5)h_{\a_4}(z^6)h_{\a_5}(z^4)h_{\a_6}(z^2),&z^3\not=1 \\
 \\
\end{array}
	\quad\quad\quad\text{for}\quad
p=2
$$
  $$
  \begin{cases}
   \begin{array}{lc}
 h_{\a_1}(-1)h_{\a_4}(-1)h_{\a_6}(-1)\quad\text{}\quad& \\
 \\
 h(z)=h_{\a_1}(z^4)h_{\a_2}(z^3)h_{\a_3}(z^5)h_{\a_4}(z^6)h_{\a_5}(z^4)h_{\a_6}(z^2),&z\not=1 \\
\end{array}
\end{cases}
	\quad\quad\quad\text{for}\quad
p=3
$$
    \end{proposition}
  \pf Let $x$ be a semisimple element of $G$, and assume $\O=\O_x$ is spherical. Without loss of generality $C(x)=L_J$, a pseudo-Levi subgroup of $G$. There exists a semisimple element $\tilde x$ in $G_{\C}$ such that $C(\tilde x)$ is $L_J$ in $G_\C$. By Lemma \ref{good}, it follows that $\O_{\tilde x}$ is a spherical semisimple conjugacy class in $G_{\C}$, and therefore, from the classification of semisimple spherical conjugacy classes in zero (or good odd) characteristic (\cite[Table 2]{CCC}, \cite[Theorem 3.6]{gio3}), it follows that
$L_J$ is of type
$A_1A_5$ or $D_5T_1$. 

Let $p=2$. Then $Z(A_1A_5)=Z(G)$ (of order 3), so that we are left with $h(z)$, for $z^3\not=1$.

Let $p=3$. Then $Z(G)=1$, and we conclude by Proposition \ref{sferiche E_6}\cvd

\noindent
We have established the information presented in Table \ref{e62} and \ref {e63}, where $w_\O=w_0\wJ$.

\begin{center}
\vskip-20pt
$$
\begin{array}{|c||c|c|c|c|}
\hline
\O  &J& w_\O & H&\dim \O \\ 
\hline
\hline 
\begin{array}{c}
h(z)=h_{\a_1}(z^4)h_{\a_2}(z^3)h_{\a_3}(z^5)h_{\a_4}(z^6)h_{\a_5}(z^4)h_{\a_6}(z^2)\\
z^3\not=1
\end{array}
& \{3,4,5\}& s_{\b_1} s_{\b_2}& D_5T_1&32
\\
\hline
\end{array}
$$
\end{center}
\begin{center} Table \totable\label{e62}:  Spherical semisimple classes in $E_6$,
$p=2$
\end{center}

\begin{center}
\vskip-20pt
$$
\begin{array}{|c||c|c|c|c|}
\hline
\O  &J& w_\O & H&\dim \O \\ 
\hline
\hline 
\begin{array}{c}
h(z)=h_{\a_1}(z^4)h_{\a_2}(z^3)h_{\a_3}(z^5)h_{\a_4}(z^6)h_{\a_5}(z^4)h_{\a_6}(z^2)\\
z\not=1
\end{array}
& \{3,4,5\}& s_{\b_1} s_{\b_2}& D_5T_1&32
\\
\hline
\begin{array}{c}
 h_{\a_1}(-1)h_{\a_4}(-1)h_{\a_6}(-1)\sim h_{\a_1}(-1) \end{array}
& \vuoto\ & w_0& A_1A_5&40
\\
\hline
\end{array}
$$
\end{center}
\begin{center} Table \totable\label{e63}:  Spherical semisimple classes in $E_6$,
$p=3$
\end{center}

\begin{proposition}\label{non-sferiche E_6} Let $\O$ be the conjugacy class of a mixed element $g$ in $E_6$, $p=2$ or $3$. Then $\O$ is not spherical. 
\end{proposition}
\pf 
Let $g=su$, Jordan decomposition. Assume that $\O$ is spherical. Then both $\O_s$ and $\O_u$ are spherical, and therefore $C(s)$ is of type $A_1A_5$ or $D_5T_1$.
A dimensional argument rules out all the possibilities except the following:

$C(s)$ is of type $A_1A_5$ and $u$ is a non-identical unipotent element in the component $A_1$ of $C(s)$ (hence $p=3$). Therefore, without loss of generality we may assume $g=h_{\a_1}(-1)x_{\a_1}(1)$, which is a regular element of the standard Levi subgroup $L_J$, for $J=\{1,2\}$. By Lemma \ref {GL(3)}, 
there is an infinite family ${\cal F}=\{g_m\mid m\in k^\ast\}\subset L_J$ such that
the $g_m$'s are all $L_J$-conjugate (hence $G$-conjugate) to $g$, and pairwise not $B_J$-conjugate. By Lemma \ref{estensione} the family $(B.g_m)_{m\in k^\ast}$ is an infinite family of (distinct) $B$-orbits, a contradiction. Hence $\O$ is not spherical.\cvd

\begin{theorem}\label{E_6}
Let $G=E_6$, $p=2$ or $3$. The spherical
  classes are either semisimple or unipotent, up to a central element if $p=2$.  The semisimple classes are represented in Table \ref{e62}, \ref{e63}, the unipotent classes are represented in \cite[Table 6, 7]{mauro-cattiva}.\cvd
\end{theorem}

\subsection{Type $E_7$.} 

Here $Z(G)=\<\t$, where $\t={h_{\a_2}(-1)h_{\a_5}(-1)h_{\a_7}(-1)}$. We put
$$
\begin{array}{lll}
\b_1 = (2,2,3,4,3,2,1),\
\b_2 =(0,1,1,2,2,2,1),\
\b_3 =(0,1,1,2,1,0,0),\\
\b_4 = \a_7,\quad
\b_5=\a_5,\quad
\b_6=\a_3,\quad
\b_7=\a_2
\end{array}
$$

\begin{proposition}\label{sferiche E_7}
 Let $x$ be an element of $E_7$, any characteristic, $\O$ its conjugacy class. If one of the following holds:
\begin{itemize}
\item[(a)] $x=h_{\a_2}(-\zeta)h_{\a_5}(\zeta)h_{\a_6}(-1)h_{\a_7}(-\zeta)$ \quad  $\zeta^2=-1$, \ $p\not=2$;
\item[(b)]  $x=h_{\a_1}(-1)$,\ $p\not=2$;
\item[(c)] $x=h(z)=h_{\a_1}(z^2)h_{\a_2}(z^3)h_{\a_3}(z^4)h_{\a_4}(z^6)h_{\a_5}(z^5)h_{\a_6}(z^4)h_{\a_7}(z^3)$ \quad for $z\not=\pm1$.
\end{itemize}
Then $\O\cap Bw_\O B\cap B^-$ is non-empty, $\dim\O=\ell(w_\O)+rk(1-w_\O)$, and $\O$ is spherical.
\end{proposition}
\pf
\noindent $(a)$\quad
Let $Y$ be the set of elements $y$ of order 4 of $T$ such that $y^2=\tau$. Then $Y$ is the disjoint union of 2 $W$-classes $Y_1$, $Y_2$: $C(y)$ is of type $A_7$ if $y\in Y_1$, of type $E_6 T_1$ if $y\in Y_2$. A representative for $Y_1$ is $h_{\a_2}(-\zeta)h_{\a_5}(\zeta)h_{\a_6}(-1)h_{\a_7}(-\zeta)$ where $\zeta$ is a square root of $-1$.

Let 
$$
y=n_{\b_1}\cdots n_{\b_7}\in w_0B
$$
$w=s_{\b_1}\cdots s_{\b_7}=w_0$.
Then $y^2=h_{\b_1}(-1)\cdots h_{\b_7}(-1)=\t$, and $\dim\O_y\geq \dim B$ by Proposition \ref{metodo}. 
Since $C(x)$ is of type $A_7$, we conclude that $x\sim y$, so that $\O\cap Bw_0 B$ is non-empty, $\dim\O=\ell(w_0)+rk(1-w_0)$ and $\O$ is spherical. As above,  $\O\cap Bw B\cap B^-\not=\emptyset$.

\noindent
$(b)$\quad
The group $G$ has 2 classes of non-central involutions: $\O_{h_{\b_1}(-1)}$ and $\O_{h_{\b_1}(-1)\tau}$: in fact there are $127$ involutions in $T$, and $\tau$ is central. The remaining 126 fall in 2 classes: $\{h_\a(-1)\mid \a\in \Phi^+\}$ and $\{h_\a(-1)\tau\mid \a\in \Phi^+\}$.
Let
$$
y=n_{\b_1}n_{\b_2}n_{\b_3}n_{\a_3}\in wB
$$
where $w=s_{\b_1}s_{\b_2}s_{\b_3}s_{\a_3}$.
Then $y^2=h_{\b_1}(-1)h_{\b_2}(-1)h_{\b_3}(-1)h_{\a_3}(-1)=1$, so that $y$ is a (non-central) involution. 
We conclude that $x\sim y$ or $x\sim y\tau$, so that (in either cases) $\O\cap Bw B$ is non-empty, $\dim\O=\ell(w)+rk(1-w)$ and $\O$ is spherical. As above,  $\O\cap Bw B\cap B^-\not=\emptyset$ (in fact we have $n_\a\sim h_{\a}(\zeta)$ already in $\<{X_\a,X_{-\a}}$ for every root $\a$, hence $n_{\b_1}n_{\b_2}n_{\b_3}n_{\a_3}\sim h_{\b_1}(\zeta)h_{\b_2}(\zeta)h_{\b_3}(\zeta)h_{\a_3}(\zeta)=h_{\gamma}(-1)$, where $\gamma=\beta_1-\a_1$. Therefore $x\sim y$).
\bigskip

\noindent 
$(c)$\quad (any characteristic)
We have $C(x)$ of type $E_6T_1$. Let 
$$
g=n_{\b_1}n_{\b_2}n_{\a_7}x_{\b_1}(1)x_{\b_2}(1)x_{\a_7}(1)
$$
Then
$$
gxg^{-1}\in \O\cap Bw B\cap B^-
$$
with 
$w=s_{\b_1}s_{\b_2}s_{\a_7}$,
$\dim\O=\ell(w)+rk(1-w)$. 
We conclude by Proposition \ref{metodo}.\cvd

  \begin{proposition}\label{semisemplici E_7}
Let $G=E_7$. The spherical semisimple
  classes are represented by
  $$
   \begin{array}{lc}
  \\
 h(z)=h_{\a_1}(z^2)h_{\a_2}(z^3)h_{\a_3}(z^4)h_{\a_4}(z^6)h_{\a_5}(z^5)h_{\a_6}(z^4)h_{\a_7}(z^3),&z\not=1 \\
 \\
\end{array}
	\ \quad\quad\quad\text{for}\quad
p=2
$$
  $$
  \begin{cases}
   \begin{array}{ll}
h_{\a_2}(-\zeta)h_{\a_5}(\zeta)h_{\a_6}(-1)h_{\a_7}(-\zeta)&\zeta^2 =-1,\\
\\
h_{\a_1}(-1), \quad h_{\a_1}(-1)\tau&\\
 \\
 h(z)=h_{\a_1}(z^2)h_{\a_2}(z^3)h_{\a_3}(z^4)h_{\a_4}(z^6)h_{\a_5}(z^5)h_{\a_6}(z^4)h_{\a_7}(z^3)&z\not=\pm1 \\
\end{array}
\end{cases}
	\text{for}\quad
p=3
$$
    \end{proposition}
  \pf Let $x$ be a semisimple element of $G$, and assume $\O=\O_x$ is spherical. Without loss of generality $C(x)=L_J$, a pseudo-Levi subgroup of $G$. There exists a semisimple element $\tilde x$ in $G_{\C}$ such that $C(\tilde x)$ is $L_J$ in $G_\C$. By Lemma \ref{good}, it follows that $\O_{\tilde x}$ is a spherical semisimple conjugacy class in $G_{\C}$, and therefore, from the classification of semisimple spherical conjugacy classes in zero (or good odd) characteristic (\cite[Table 2]{CCC}, \cite[Theorem 3.7]{gio3}), it follows that
$L_J$ is of type
$E_6T_1$, $D_6A_1$ or $A_7$. 

Let $p=2$. Then $Z(G)=Z(D_6A_1)=Z(A_7)=1$, so that we are left with $h(z)$, for $z\not=1$.

Let $p=3$: we conclude by Proposition \ref{sferiche E_7}\cvd

\noindent
We have established the information presented in Table \ref{e62} and \ref {e63}, where $w_\O=w_0\wJ$.

\begin{center}
\vskip-20pt
$$
\begin{array}{|c||c|c|c|c|}
\hline
\O  &J& w_\O & H&\dim \O \\ 
\hline
\hline 
\begin{array}{c}
h(z)=h_{\a_1}(z^2)h_{\a_2}(z^3)h_{\a_3}(z^4)h_{\a_4}(z^6)h_{\a_5}(z^5)h_{\a_6}(z^4)h_{\a_7}(z^3)\\
z\not=1
\end{array}
& \{2,3,4,5\}& s_{\b_1} s_{\b_2}s_{\a_7}& E_6T_1&54
\\
\hline 
\end{array}
$$
\end{center}
\begin{center} Table \totable\label{e72}:  Spherical semisimple classes in $E_7$,
$p=2$
\end{center}

\begin{center}
\vskip-20pt
$$
\begin{array}{|c||c|c|c|c|}
\hline
\O  &J& w_\O & H&\dim \O \\ 
\hline
\hline 
\begin{array}{c}
h(z)=h_{\a_1}(z^2)h_{\a_2}(z^3)h_{\a_3}(z^4)h_{\a_4}(z^6)h_{\a_5}(z^5)h_{\a_6}(z^4)h_{\a_7}(z^3)\\
z\not=\pm1
\end{array}
& \{2,3,4,5\}& s_{\b_1} s_{\b_2}s_{\a_7}& E_6T_1&54
\\
\hline 
\begin{array}{c}
h_{\a_1}(-1),  h_{\a_1}(-1)\tau
\end{array}
& \{2,5,7\}& s_{\b_1} s_{\b_2}s_{\b_3}s_{\a_3}& D_6A_1&64
\\
\hline
\begin{array}{c}
h_{\a_2}(-\zeta)h_{\a_5}(\zeta)h_{\a_6}(-1)h_{\a_7}(-\zeta)\\
\zeta^2=-1
\end{array}
& \vuoto\ & w_0& A_7&70
\\
\hline
\end{array}
$$
\end{center}
\begin{center} Table \totable\label{e73}:  Spherical semisimple classes in $E_7$,
$p=3$
\end{center}

\begin{proposition}\label{non-sferiche E_7} Let $\O$ be the conjugacy class of a mixed element $g$ in $E_7$, $p=2$ or $3$. Then $\O$ is not spherical. 
\end{proposition}
\pf 
Let $g=su$, Jordan decomposition. Assume that $\O$ is spherical. Then both $\O_s$ and $\O_u$ are spherical, and therefore $C(s)$ is of type $E_6T_1$, $D_6A_1$ or $A_7$.
A dimensional argument rules out all the possibilities except the following:

$C(s)$ is of type $D_6A_1$ and $u$ is a non-identical unipotent element in the component $A_1$ of $C(s)$ (hence $p=3$). Therefore, without loss of generality we may assume $g=h_{\a_7}(-1)x_{\a_7}(1)$, which is a regular element of the standard Levi subgroup $L_J$, for $J=\{6,7\}$. By Lemma \ref {GL(3)}, 
there is an infinite family ${\cal F}=\{g_m\mid m\in k^\ast\}\subset L_J$ such that
the $g_m$'s are all $L_J$-conjugate (hence $G$-conjugate) to $g$, and pairwise not $B_J$-conjugate. By Lemma \ref{estensione} the family $(B.g_m)_{m\in k^\ast}$ is an infinite family of (distinct) $B$-orbits, a contradiction. Hence $\O$ is not spherical.\cvd

\begin{theorem}\label{E_7}
Let $G=E_7$, $p=2$ or $3$. The spherical
  classes are either semisimple or unipotent, up to a central element if $p=3$ .  The semisimple classes are represented in Table \ref{e72}, \ref{e73}, the unipotent classes are represented in \cite[Table 8, 9]{mauro-cattiva}.\cvd
\end{theorem}

\subsection{Type $E_8$.}  

We put
$$
\begin{array}{lll}
&\b_1 =   (2,3,4,6,5,4,3,2),\ 
\b_2=(2,2,3,4,3,2,1,0),\
\b_3=(0,1,1,2,2,2,1,0),\\
&\b_4=(0,1,1,2,1,0,0,0),\
\b_5=\a_7, \
\b_6=\a_5,\
\b_7=\a_3,\
\b_8=\a_2
\end{array}
$$

\begin{proposition}\label{sferiche E_8}
 Let $x$ be an element of $E_8$, $p\not=2$, $\O$ its conjugacy class. If one of the following holds:
\begin{itemize}
\item[(a)] $x=h_{\a_2}(-1)h_{\a_3}(-1)$;
\item[(b)]  $x=h_{\a_2}(-1)h_{\a_5}(-1)h_{\a_7}(-1)\sim h_{\a_8}(-1)$.
\end{itemize}
Then $\O\cap Bw_\O B\cap B^-$ is non-empty, $\dim\O=\ell(w_\O)+rk(1-w_\O)$, and $\O$ is spherical.
\end{proposition}
\pf The group $E_8$, for $p\not=2$, has 2 classes of involutions.

\noindent $(a)$\quad
Let 
$$
y=n_{\b_1}\cdots n_{\b_8}\in w_0B
$$
$w=s_{\b_1}\cdots s_{\b_8}=w_0$.
Then $y^2=h_{\b_1}(-1)\cdots h_{\b_8}(-1)=1$, and $\dim\O_y\geq \dim B$ by Proposition \ref{metodo}. 
Since $C(x)$ is of type $D_8$, we conclude that $x\sim y$, so that $\O\cap Bw_0 B$ is non-empty, $\dim\O=\ell(w_0)+rk(1-w_0)$ and $\O$ is spherical. As above,  $\O\cap Bw B\cap B^-\not=\emptyset$.

\bigskip

\noindent 
$(b)$\quad 
Let $x=h_{\a_8}(-1)$, so that 
$C(x)$ is of type $E_7A_1$. Let 
$$
g=n_{\b_1}n_{\b_2}n_{\b_3}n_{\a_7}x_{\b_1}(1)x_{\b_2}(1)x_{\b_3}(1)x_{\a_7}(1)
$$
Then
$$
gxg^{-1}\in \O\cap Bw B\cap B^-
$$
with 
$w=s_{\b_1}s_{\b_2}s_{\b_3}s_{\a_7}$,
$\dim\O=\ell(w)+rk(1-w)$. 
We conclude by Proposition \ref{metodo}.\cvd

  \begin{proposition}\label{semisemplici E_8}
Let $G=E_8$. The spherical (non-trivial) semisimple
  classes are represented by
  $$
   \begin{array}{lc}
  \\
 \text{none}& \\
 \\
\end{array}
\quad	\quad \quad\quad\quad\text{for}\quad
p=2
$$
  $$
  \begin{cases}
   \begin{array}{ll}
h_{\a_2}(-1)h_{\a_3}(-1)&\\
\\
h_{\a_8}(-1)&\\
\end{array}
\end{cases}
	\text{for}\quad
p=3\quad \text{or}\quad 5
$$
    \end{proposition}
  \pf Let $x$ be a semisimple element of $G$, and assume $\O=\O_x$ is spherical. Without loss of generality $C(x)=L_J$, a pseudo-Levi subgroup of $G$. There exists a semisimple element $\tilde x$ in $G_{\C}$ such that $C(\tilde x)$ is $L_J$ in $G_\C$. By Lemma \ref{good}, it follows that $\O_{\tilde x}$ is a spherical semisimple conjugacy class in $G_{\C}$, and therefore, from the classification of semisimple spherical conjugacy classes in zero (or good odd) characteristic (\cite[Table 2]{CCC}, \cite[Theorem 3.8]{gio3}), it follows that
$L_J$ is of type
$E_7A_1$ or $D_8$. 

Let $p=2$. Then $Z(E_7A_1)=Z(D_8)=1$, so there are no non-trivial spherical semisimple classes.

Let $p=3$ or $5$: we conclude by Proposition \ref{sferiche E_8}.\cvd

\noindent
We have established the information presented in Table \ref{e62} and \ref {e63}, where $w_\O=w_0\wJ$.

\begin{center}
\vskip-20pt
$$
\begin{array}{|c||c|c|c|c|}
\hline
\O  &J& w_\O & H&\dim \O \\ 
\hline
\hline 
\begin{array}{c}
h_{\a_8}(-1)
\end{array}
& \{2,3,4,5\}& s_{\b_1} s_{\b_2}s_{\b_3}s_{\a_7}& E_7A_1&112
\\
\hline
\begin{array}{c}
h_{\a_2}(-1)h_{\a_3}(-1)
\end{array}
& \vuoto\ & w_0& D_8&128
\\
\hline
\end{array}
$$
\end{center}
\begin{center} Table \totable\label{e835}:  Spherical semisimple classes in $E_8$,
$p=3$ or $5$
\end{center}

\begin{proposition}\label{non-sferiche E_8} Let $\O$ be the conjugacy class of a mixed element $g$ in $E_8$, $p=2$, $3$ or $5$. Then $\O$ is not spherical. 
\end{proposition}
\pf 
Let $g=su$, Jordan decomposition. Assume that $\O$ is spherical. Then both $\O_s$ and $\O_u$ are spherical, and therefore $C(s)$ is of type $E_7A_1$ or $A_7$.
A dimensional argument rules out all the possibilities except the following:

$C(s)$ is of type $E_7A_1$ and $u$ is a non-identical unipotent element in the component $A_1$ of $C(s)$ (hence $p=3$ or $5$). Therefore, without loss of generality we may assume $g=h_{\a_8}(-1)x_{\a_8}(1)$, which is a regular element of the standard Levi subgroup $L_J$, for $J=\{7,8\}$. By Lemma \ref {GL(3)}, 
there is an infinite family ${\cal F}=\{g_m\mid m\in k^\ast\}\subset L_J$ such that
the $g_m$'s are all $L_J$-conjugate (hence $G$-conjugate) to $g$, and pairwise not $B_J$-conjugate. By Lemma \ref{estensione} the family $(B.g_m)_{m\in k^\ast}$ is an infinite family of (distinct) $B$-orbits, a contradiction. Hence $\O$ is not spherical.\cvd

\begin{theorem}\label{E_8}
Let $G=E_8$, $p=2$, $3$ or $5$. The spherical
  classes are either semisimple or unipotent.  The semisimple classes are represented in Table \ref{e835}, the unipotent classes are represented in \cite[Table 10, 11]{mauro-cattiva}.\cvd
\end{theorem}

\subsection{Type $F_4$.}  

We put
$$
\begin{array}{ll}
\b_1 =(2,3,4,2),&\b_2=(0,1,2,2),\\
\b_3=(0,1,2,0),&\b_4=(0,1,0,0)
\end{array}
$$
Moreover $\gamma_1$ is the highest short root $(1,2,3,2)$.

\begin{proposition}\label{sferiche F_4}
 Let $x$ be an element of $F_4$, $p\not=2$, $\O$ its conjugacy class. If one of the following holds:
\begin{itemize}
\item[(a)] $x=h_{\a_2}(-1)h_{\a_4}(-1)\sim h_{\a_1}(-1)$;
\item[(b)]  $x=h_{\a_4}(-1)$.
\end{itemize}
Then $\O\cap Bw_\O B\cap B^-$ is non-empty, $\dim\O=\ell(w_\O)+rk(1-w_\O)$, and $\O$ is spherical.
\end{proposition}
\pf The group $F_4$, for $p\not=2$, has 2 classes of involutions.

\noindent $(a)$\quad
Let 
$$
y=n_{\b_1}\cdots n_{\b_4}\in w_0B
$$
$w=s_{\b_1}\cdots s_{\b_4}=w_0$.
Then $y^2=h_{\b_1}(-1)\cdots h_{\b_4}(-1)=1$, and $\dim\O_y\geq \dim B$ by Proposition \ref{metodo}. 
Since $C(x)$ is of type $A_1C_3$, we conclude that $x\sim y$, so that $\O\cap Bw_0 B$ is non-empty, $\dim\O=\ell(w_0)+rk(1-w_0)$ and $\O$ is spherical. As above,  $\O\cap Bw B\cap B^-\not=\emptyset$.

\bigskip

\noindent 
$(b)$\quad 
We have $C(x)$ of type $B_4$. Let 
$$
g=n_{\gamma_1}x_{\gamma_1}(1)
$$
Then
$$
gxg^{-1}\in \O\cap Bw B\cap B^-
$$
with 
$w=s_{\gamma_1}$,
$\dim\O=\ell(w)+rk(1-w)$. 
We conclude by Proposition \ref{metodo}.\cvd

  \begin{proposition}\label{semisemplici F_4}
Let $G=F_4$. The spherical (non-trivial) semisimple
  classes are represented by
  $$
   \begin{array}{lc}
  \\
 \text{none}& \\
 \\
\end{array}
\quad	\quad \quad\quad\quad\text{for}\quad
p=2
$$
  $$
  \begin{cases}
   \begin{array}{ll}
h_{\a_1}(-1)&\\
\\
h_{\a_4}(-1)&\\
\end{array}
\end{cases}
	\text{for}\quad
p=3
$$
    \end{proposition}
  \pf Let $x$ be a semisimple element of $G$, and assume $\O=\O_x$ is spherical. Without loss of generality $C(x)=L_J$, a pseudo-Levi subgroup of $G$. There exists a semisimple element $\tilde x$ in $G_{\C}$ such that $C(\tilde x)$ is $L_J$ in $G_\C$. By Lemma \ref{good}, it follows that $\O_{\tilde x}$ is a spherical semisimple conjugacy class in $G_{\C}$, and therefore, from the classification of semisimple spherical conjugacy classes in zero (or good odd) characteristic (\cite[Table 2]{CCC}, \cite[Theorem 3.9]{gio3}), it follows that
$L_J$ is of type
$C_3A_1$ or $B_4$. 

Let $p=2$. Then $Z(C_3A_1)=Z(B_4)=1$, so there are no non-trivial spherical semisimple classes.

Let $p=3$: we conclude by Proposition \ref{sferiche F_4}.\cvd

\noindent
We have established the information presented in Table \ref{f43}, where $w_\O=w_0\wJ$.

\begin{center}
\vskip-20pt
$$
\begin{array}{|c||c|c|c|c|}
\hline
\O  &J& w_\O & H&\dim \O \\ 
\hline
\hline 
\begin{array}{c}
h_{\a_4}(-1)
\end{array}
& \{1,2,3\}& s_{\gamma_1} & B_4&16
\\
\hline
\begin{array}{c}
h_{\a_1}(-1)
\end{array}
& \vuoto\ & w_0& C_3A_1&28
\\
\hline
\end{array}
$$
\end{center}
\begin{center} Table \totable\label{f43}:  Spherical semisimple classes in $F_4$,
$p=3$
\end{center}

\begin{proposition}\label{mista sferica F_4}
 Let $x=h_{\a_4}(-1)x_{\a_4}(1)$ in $F_4$, $p\not=2$, $\O$ its conjugacy class. Then $\O\cap Bw_\O B\cap B^-$ is non-empty, $\dim\O=\ell(w_\O)+rk(1-w_\O)$, and $\O$ is spherical.
\end{proposition}
\pf 
This is the mixed class in $F_4$ which is spherical in zero or good, odd characteristic.
We can deal with this class with the same method used in the proof of \cite[Theorem 23]{CCC}, and corrected in the proof of \cite[Theorem 3.9]{gio3}, to show that $Bw_0B\cap\O\not=\emptyset$, so that $w_{\O}=w_0$, $\dim\O=\ell(w_0)+rk (1-w_0)$, and $\O$ is spherical. The (correct) argument at the end of the proof of \cite[Theorem 23]{CCC}, shows that $\O\cap Bw_0B\cap B^-\not=\emptyset$.\cvd

 \begin{proposition}\label{miste F_4}
Let $G=F_4$. The spherical mixed
  classes are represented by
  $$
   \begin{array}{lc}
 \text{none}& 
 \text{for}\quad
p=2
\\
 \\
 h_{\a_4}(-1)x_{\a_4}(1)&\text{for}\quad
p=3
\end{array}
$$
    \end{proposition}
  \pf
  Let $g=su$, Jordan decomposition, be a mixed element. Assume that $\O=\O_g$ is spherical. Then both $\O_s$ and $\O_u$ are spherical, and therefore $C(s)$ is of type $C_3A_1$ or $B_4$.
A dimensional argument rules out all the possibilities except the following:

$C(s)$ is of type $B_4$ and $u$ is in the minimal unipotent class of $C(s)$ (hence $p=3$). 
Therefore, without loss of generality we may assume $g=h_{\a_4}(-1)x_{\a_4}(1)$. We conclude by \ref{mista sferica F_4}.
\cvd

\begin{center}
\vskip-20pt
$$
\begin{array}{|c||c|c|c|c|}
\hline
\O  &J& w_\O & \dim \O \\ 
\hline
\hline 
\begin{array}{c}
h_{\a_4}(-1)x_{\a_4}(1)\\
\end{array}
& \vuoto& w_0& 28
\\
\hline 
\end{array}
$$
\end{center}
\begin{center} Table \totable\label{f4m}:  Spherical mixed classes in $F_4$,
$p=3$
\end{center}

\begin{theorem}\label{F_4}
Let $G=F_4$, $p=2$ or $3$. If $p=2$ the spherical classes are unipotent, and are represented in \cite[Table 13]{mauro-cattiva}. If $p=3$, 
the spherical semisimple classes are represented in Table \ref{f43}, the spherical unipotent classes are represented in \cite[Table 12]{mauro-cattiva}, the spherical mixed classes are represented in Table \ref{f4m}.\cvd
\end{theorem}

\medskip

\subsection{Type $G_2$.}

We put 
$
\b_1 =(3,2),\
\b_2=\a_1 
$.
Moreover $\gamma_1$ is the highest short root $(2,1)$.

\medskip

\begin{proposition}\label{sferiche G_2}
 Let $x$ be an element of $G_2$, any characteristic, $\O$ its conjugacy class. If one of the following holds:
\begin{itemize}
\item[(a)] $x=h_{\a_1}(-1)$, \ $p\not=2$;
\item[(b)]  $x=h_{\a_1}(\zeta)$,\ $\zeta$ a primitive $3$-rd root of 1, \ $p\not=3$;
\end{itemize}
Then $\O\cap Bw_\O B\cap B^-$ is non-empty, $\dim\O=\ell(w_\O)+rk(1-w_\O)$, and $\O$ is spherical.
\end{proposition}
\pf 

\noindent $(a)$\quad
$p\not=2$. The group $G_2$ has 1 class of involutions.
Let 
$$
y=n_{\b_1}n_{\b_2}\in w_0B
$$
$w=s_{\b_1}s_{\b_2}=w_0$.
Then $y^2=h_{\b_1}(-1) h_{\b_2}(-1)=1$, and $\dim\O_y\geq \dim B$ by Proposition \ref{metodo}. 
We conclude that $x\sim y$, so that $\O\cap Bw_0 B$ is non-empty, $\dim\O=\ell(w_0)+rk(1-w_0)$ and $\O$ is spherical. As above,  $\O\cap Bw B\cap B^-\not=\emptyset$.

\bigskip

\noindent 
$(b)$\quad 
$p\not=3$.
Let 
$$
g=n_{\g_1}x_{\g_1}(1)
$$
Then
$$
gxg^{-1}\in \O\cap Bw B\cap B^-
$$
with 
$w=s_{\g_1}$,
$\dim\O=\ell(w)+rk(1-w)$. 
We conclude by Proposition \ref{metodo}.\cvd

  \begin{proposition}\label{semisemplici G_2}
Let $G=G_2$. The spherical semisimple classes are represented by
$$
\begin{array}{lc}
\\
\ \ h_{\a_1}(-1)& \text{for}\quad p=3
\\
\\
\begin{array}{lc}
h_{\a_1}(\zeta)\\ \zeta\ \text{a primitive $3$-rd root of}\ \ 1 
\end{array}&
\text{for}\quad
p=2
\end{array}
$$
    \end{proposition}
  \pf Let $x$ be a semisimple element of $G$, and assume $\O=\O_x$ is spherical. Without loss of generality $C(x)=L_J$, a pseudo-Levi subgroup of $G$. There exists a semisimple element $\tilde x$ in $G_{\C}$ such that $C(\tilde x)$ is $L_J$ in $G_\C$. By Lemma \ref{good}, it follows that $\O_{\tilde x}$ is a spherical semisimple conjugacy class in $G_{\C}$, and therefore, from the classification of semisimple spherical conjugacy classes in zero (or good odd) characteristic (\cite[Table 2]{CCC}, \cite[Theorem 3.1]{gio3}), it follows that
$L_J$ is of type
$A_1\tilde A_1$ or $A_2$. If $p=2$, then $Z(A_1\tilde A_1)=1$, if $p=3$, then $Z(A_2)=1$: we conclude by \ref{sferiche G_2}.\cvd

We have established the information presented in Table \ref{g22} and \ref {g23}, where $w_\O=w_0\wJ$.

\begin{center}
\vskip-20pt
$$
\begin{array}{|c||c|c|c|c|}
\hline
\O  &J& w_\O & H&\dim \O \\ 
\hline
\hline 
\begin{array}{c}
h_{\a_1}(\zeta)
\\
\zeta\ \text{a primitive 3-rd root of 1}\end{array}
& \{2\}& s_{\gamma_1} & A_2&6
\\
\hline
\end{array}
$$
\end{center}
\begin{center} Table \totable\label{g22}:  Spherical semisimple classes in $G_2$,
$p=2$
\end{center}

\begin{center}
\vskip-20pt
$$
\begin{array}{|c||c|c|c|c|}
\hline
\O  &J& w_\O & H&\dim \O \\ 
\hline
\hline 
\begin{array}{c}
\\
h_{\a_1}(-1)\\
\\
\end{array}
& \vuoto\ & w_0& A_1\tilde A_1&8
\\
\hline
\end{array}
$$
\end{center}
\begin{center} Table \totable\label{g23}:  Spherical semisimple classes in $G_2$,
$p=3$
\end{center}

\begin{theorem}\label{G2}
Let $G=G_2$, $p=2$, $3$. The spherical
  classes are either semisimple or unipotent.  The semisimple classes are represented in Table \ref{g22}, \ref{g23}, the unipotent classes are represented in \cite[Table 14, 15]{mauro-cattiva}.
  \end{theorem}
  \pf By the above discussion, we are left to show that no mixed class is spherical. Let $g=su$, Jordan decomposition. Assume that $\O_g$ is spherical. Then both $\O_s$ and $\O_u$ are spherical, and therefore $C(s)$ is of type $A_1\tilde A_1$ or $A_2$.
A dimensional argument rules out all the possibilities.
  \cvd

\section{Final remarks}

Once achieved the classification of spherical conjugacy classes in all characteristics, and proved that for every spherical conjugacy class $\O$ we have $\dim\O=\ell(w_\O)+rk (1-w_\O)$ one can extend results obtained in \cite{gio1}, \cite{gio2}, \cite{CCC}, \cite{lu}. In \cite[Theorem 25]{CCC} we established the characterization of spherical conjugacy classes in terms of the dimension formula: a conjugacy class $\O$ in $G$ is spherical if and only if $\dim{\cal O}=\ell(w_\O)+rk (1-w_\O)$. This was obtained over the complex numbers, and the same proof works over any algebraically closed field of characteristic zero. Then the same characterization was given in zero, or good odd characteristic in \cite[Theorem 4.4]{gio1}, without the classification of spherical conjugacy classes. Lu gave a very neat proof of the dimension formula (even for twisted conjugacy classes) in  \cite[Theorem 1.1]{lu} in characteristic zero. From the results obtained in the previous section, we may state:

\begin{theorem}\label{finale} Let $\O$ be a conjugacy class of a simple algebraic group, any characteristic. 
The following are equivalent:
\begin{itemize}
\item[(a)] $\O$ is spherical;
\item[(b)] there exists $w\in W$ such that $\O\cap BwB\not=\emptyset$ and $\dim\O\leq\ell(w)+\rk (1-w)$;
\item[(c)]  $\dim{\cal O}=\ell(w_\O)+rk (1-w_\O)$.\cvd
\end{itemize}
\end{theorem}

\begin{corollary}\label{} Let $\O$ be a spherical class of $G$. Then $\dim{\cal O}\leq\ell(w_0)+rk (1-w_0)$.
\end{corollary}
\pf We have $\dim{\cal O}=\ell(w_\O)+rk (1-w_\O)$, and 
$\ell(w)+rk (1-w)\leq \ell(w_0)+rk (1-w_0)$ for every $w\in W$ (cf. \cite[Remark 4.14]{gio1}).\cvd
\bigskip

\begin{proposition}\label{centralizzanti in wB} Let $\O$ be a spherical conjugacy class, $w=w_\O=w_0 \wJ$. Then $(T^{w})^\circ \leq C_T(x)\leq T^{w}$, $C_U(x)=U_{\wJ}$ and $(T^{w})^\circ U_{\wJ}\leq C_B(x)\leq T^{w} U_{\wJ}$  for every $x\in w B$. 
\end{proposition}
\pf We choose a representative $\dot w$ of $w$ in $N$ such that $x=\dot w u$, $u\in U$. Let $b\in C_B(x)$, $b=t u_1 u_2$, where $t\in T$, $u_1\in U_{w}$, $u_2\in U_{\wJ}$.  From the Bruhat decomposition, we get $u_1=1$, $t\in T^{w}$, so that $C_B(x)\leq T^{w} U_{\wJ}$. But $\dim \O=\ell(w)+rk (1-w)$ implies $\dim C_B(x)=n-rk(1-w)+\ell(\wJ)=\dim T^{w}U_{\wJ}$. Hence $(C_B(x))^\circ=(T^{w})^\circ U_{\wJ}$ and $C_U(x)=U_{\wJ}$.

Now assume $b\in C_{TU_{w}}(x)$, $b=t u_1$, where $t\in T$, $u_1\in U_{w}$. Again from the Bruhat decomposition, we get $u_1=1$, $t\in T^{w}$, so that $C_{TU_{w}}(x)=C_{T}(x)\leq T^{w}$.  We have
$B.x=TU_{w}U_{\wJ}.x=TU_{w}.x$, hence $\dim C_{TU_{w}}(x)=n-rk(1-w)=\dim T^w$. It follows that 
$(T^{w})^\circ \leq C_T(x)\leq T^{w}$.\cvd

\begin{theorem}\label{centralizzanti} Let $\O$ be a spherical conjugacy class of a simple algebraic group, $v=\O\cap Bw_\O B$ the dense $B$-orbit. Then $C_U(x)$ is connected and 
 $C_B(x)$ is a split extension of $(C_B(x))^\circ$ by an elementary abelian $2$-group for every $x\in v$.
If $p=2$, then $C_U(x)$, $C_T(x)$ and $C_B(x)$ are connected for every $x\in \O \cap Bw_\O B$.
\end{theorem}
\pf Let $w=w_\O$. We may assume $x\in wB$. Following the discussion after \cite[Corollary 3.22]{mauro-mathZ}, we have $T=(T^w)^\circ(S^w)^\circ$, where $S^w=\{t\in T\mid t^w=t^{-1}\}$. Then $T^w=(T^w)^\circ(T^w\cap T_2)$, where $T_2=\{t\in T\mid t^2=1\}$, and $C_T(x)=(T^w)^\circ C_{T_2}(x)$  by Proposition \ref{centralizzanti in wB}. There exists a subgroup $R$ of $T_2$ such that $T^w=(T^w)^\circ \times R$, hence $C_T(x)=(T^w)^\circ\times C_{R}(x)$. In particular, $C_B(x)=((T^w)^\circ\times C_{R}(x))U_{\wJ}=(C_B(x))^\circ\ C_{R}(x)$. If $p=2$, then $T_2=\{1\}$, $T^w=(T^w)^\circ=C_T(x)$ and $C_B(x)=(C_B(x))^\circ$.\cvd

\medskip

We recall, from Remark \ref{odd}, that there is an action of $W$ on the set ${\cal V}$ of $B$-orbits in $\O$ when $\O$ is a spherical conjugacy class and $p\not=2$. We are in the position to prove that this action is defined also for $p=2$.

\begin{corollary}\label{}  Let $\O$ be a spherical conjugacy class of a simple algebraic group, any characteristic. Then there is an action of the Weyl group $W$ on the set of $B$-orbits in $\O$, as defined in \cite{knop}.
\end{corollary}
\pf We have only to deal with $p=2$.
By \cite[Theorem 4.2 c)]{knop}, the action of $W$ is defined on the set of $B$-orbits in $\O$ as long as 
$C_U(x)$ is connected for every $x\in \O$. By Theorem \ref{centralizzanti}, $C_U(x)$ is connected for every $x$ in the dense $B$-orbit: this ensures that $C_U(x)$ is connected for every $x\in\O$, by \cite[Corollary 3.4)]{knop}.\cvd

\medskip

\def\wK{w_{\!_K}}

Once the $W$-action has been defined when $p=2$, we can extend to this case the results obtained by G. Carnovale in zero or good odd characteristic.

\begin{theorem}\label{Giovanna} Let $\O$ be a spherical conjugacy class of a simple algebraic group. If $\O\cap BwB$ is non-empty, then $w^2=1$.
\end{theorem}
\pf We may apply the same argument ued in the proof of \cite[Theorem 2.7]{gio1}, which uses both the monoid action of $M(W)$ and the group action of $W$ on the set of $B$-orbits in $\O$.\cvd

\begin{theorem}\label{invo} Let $\O$ be a conjugacy class in a simple algebraic group. If 
$$
\{w\in W\mid \O \cap BwB\not=\emptyset\}\subseteq\{w\in W\mid w^2=1\}
$$
 then $\O$ is spherical.
\end{theorem}
\pf Again, we may apply the proof of 
 \cite[Theorem 5.7]{gio2}: note that in characteristic 2 one does not need to consider groups of type $B_n$. \cvd

\medskip
Assume  $\O$ is a spherical conjugacy class of a simple algebraic group (any characteristic), $v$ the dense $B$-orbit in $\O$. Set  $P=\{g\in G\mid g.v=v\}$. Then $P$ is a parabolic subgroup of $G$ containing $B$, and therefore $P=P_K$, the standard parabolic subgroup relative to a certain subset $K$ of $\Pi$. 

\begin{theorem}\label{} Let $\O$ be a spherical conjugacy class of a simple algebraic group, any characteristic, $w=w_0\wJ$ be the unique element in $W$ such that $\O\cap BwB$ is dense in $\O$, $v=\O\cap BwB$ the dense $B$-orbit in $\O$, $P_K=\{g\in G\mid g.v=v\}$. Then $K=J$. If $x\in \O\cap wB$, then $L_J'$ and $(T^w)^\circ$ are contained in $C(x)$ and $C_B(x)^\circ = (T^w)^\circ U_{\wJ}$.
\end{theorem}
\pf We have already showed that $C_B(x)^\circ = (T^w)^\circ U_{\wJ}$ for every $x\in \O\cap wB$.
Let $S=\{i,\vartheta(i)\}$ be a $\vartheta$-orbit in $\Pi\setminus J$ consisting of 2 elements. We put $H_S=\{h_{\a_i}(z)h_{\a_{\vartheta(i)}}(z^{-1})\mid z\in k^\ast\}
$. Let ${\cal S}_1$ be the set of $\vartheta$-orbits in $\Pi\setminus J$ consisting of 2 elements.  Then, by \cite[Remark 3.10]{mauro-mathZ}, $\Delta_J\cup \{\a_i-\a_{\vartheta(i)}\}_{{\cal S}_1}$ is a basis of $\ker(1-w)$ and
\begin{equation}
(T^w)^\circ=\prod_{j\in J}H_{\a_j}\times \prod_{S\in {\cal S}_1}
H_S
\label{prodotto}
\end{equation}
We put $\Psi_J=\{\beta\in \Phi\mid w(\b)=-\b\}$. Then $\Psi_J$ is a root system in $\im(1-w)$ (\cite[Proposition 2]{springer2}), and $w_{|\im(1-w)}$ is $-1$.
If $K=C((T^w)^\circ)'$, then $K$ is semisimple with root system $\Psi_J$ and maximal torus $T\cap K=(S^w)^\circ$. Assume $x=\dot w u\in v$, with $u\in U$. Then $(T^w)^\circ\leq C(x)$ implies $x\in C((T^w)^\circ)$, moreover $\dot w \in C(T^w)$, so that $u\in K$. Let $u=\prod_{\a\in \Phi^+\cap \Psi_J}x_\a(k_\a)$ be the expression of $u$ for any fixed total ordering on $\Phi^+$. If $k_\a\not=0$, then $w(\a)=-\a$, so that in particular $u\in U_w$. Moreover, if $\beta\in \Phi_J$, then $(\a,\b)=(w\a,w\b)=(-\a,\b)$, so that $\a\perp\b$. Finally, we have $\vartheta\a=-\a$, since $w\a=-\a$ is equivalent to $\wJ\a=-w_0\a$, and $\wJ\a=\a$, since $\wJ\in W_J$, and $(\a,\a_j)=0$ for every $j\in J$.

From the fact that $U_{\wJ}\leq C(x)$, it follows that $U_{\wJ}\leq C(\dot w)$, hence $U_{\wJ}\leq C(u)$.
From Chevalley commutator formula, we deduce that $\wJ U_{\wJ}\wJ^{-1}\leq C(u)$ as well, hence $L_J'\leq C(x)$. Then we may argue as in the proof of \cite[Proposition 4.15]{gio1} to conclude that $K=J$.\cvd

\begin{remark}\label{reductive} 
{\rm 
Assume $G$ is a connected reductive algebraic group over $\K$. From the classification of spherical conjugacy classes obtained in simple algebraic groups (which is independent of the isogeny class), one gets the classification of spherical conjugacy classes in $G$. In fact, If $G=ZG_1\cdots G_r$, where $Z$ is the connected component fo the center of $G$, and $G_1,\ldots,G_r$ are the simple components of $G$, then the conjugacy class $\O$ in $G$ of $x=zx_1\cdots x_r$, with $z\in Z$, $x_i\in G_i$ for $i=1,\ldots,r$ is spherical, if and only if the conjugacy class $\O_i$ of $x_i$ in $G_i$ is spherical for every $i=1,\ldots,r$.}\cvd
\end{remark}

\begin{remark}\label{B meno} 
{\rm In order to show that a conjugacy class $\O$ is spherical, we showed that $\dim\O=\ell(w_\O)+rk(1-w_\O)$. However, in each case we even showed that $\O\cap Bw_\O B\cap B^-\not=\emptyset$. The motivation for this was the proof of the De Concini-Kac-Procesi for quantum groups at roots of one over spherical conjugacy classes (see \cite{CCC}). The fact that $\O\cap Bw_\O B\cap B^-\not=\emptyset$ for every spherical conjugacy class has been proved in characteristic zero (in \cite{CCC}). It is a general fact that if $\O$ is semisimple, then $\O\cap Bw B\not=\emptyset$ implies $\O\cap BwB\cap B^-\not=\emptyset$ for any $w\in W$ (\cite[Lemma 14]{CCC}). For unipotent classes, we showed in \cite{mauro-cattiva} that when $p=2$ then $\O\cap Bw_\O B\cap B^-\not=\emptyset$ by exhibiting explicitely an element in $\O\cap Bw_\O B\cap B^-$. The argument in  \cite[Lemma 10]{CCC}, allows to prove that $\O\cap Bw_\O B\cap B^-\not=\emptyset$ for every spherical unipotent class in good characteristic. However, it is possible to adapt the same proof to the remaining unipotent classes in bad characteristic, due to the fact that we do have the classification, and so we just make a case by case consideration. Assume $\O$ is a spherical mixed class. In all cases, apart from $F_4$, we have an explicit element in $\O\cap Bw_\O B\cap B^-$. We observed in Proposition \ref {mista sferica F_4} that the argument used in \cite{CCC} holds for every odd characteristic. We conclude that in all characteristics, if $\O$ is a spherical conjugacy class, then $\O\cap Bw_\O B\cap B^-\not=\emptyset$.
}\cvd
\end{remark}

\section{Appendix}

In this appendix we prove Lemma \ref{GL(3)} and \ref{Sp(4)}. In both cases the existence of an infinite family of $B$-orbits in the relevant conjugacy class $\O$ follows from the fact that in each case $\O$ is not spherical, by Lemma \ref{non-sferica}. However, one can easily find explicit representatives for such classes.

\begin{lemma}Let $H=GL(3)$, any characteristic, $g$ a regular element of $H$. Then there exists a subset
${\cal F}=\{x_m\mid m\in k^\ast\}$ of  $\O_g$ such that $(B(H).x_m)_{ m\in k^\ast}$ consists of pairwise distinct $B(H)$-orbits.
\end{lemma}
\pf
For $m$, $a$, $b$, $c\in k^\ast$, let
$$
x_m=x_m(a,b,c)=
\left(
\begin{array}{ccc}
 0 & 0 & \frac{a b c}{m} \\
 0 & -m & -\frac{(a+m) (b+m) (c+m)}{m} \\
 1 & 1 & a+b+c+m \\
\end{array}
\right)
=
$$
$$
=
\left(
\begin{array}{ccc}
 0 & 0 & 1 \\
 0 & -1 & 0 \\
 1 & 0 & 0 \\
\end{array}
\right)
\left(
\begin{array}{ccc}
 1 & 0 & 0 \\
 0 & m & 0 \\
 0 & 0 & \frac{a b c}{m} \\
\end{array}
\right)
\left(
\begin{array}{ccc}
 1 & 1 & a+b+c+m \\
 0 & 1 & \frac{(a+m) (b+m) (c+m)}{m^2} \\
 0 & 0 & 1 \\
\end{array}
\right)\in w_0 B
$$
We have $B.x_m\cap w_0B=T.x_m$,
$C_T(x_m)=$ scalar matrices, and
$$
S=\{
\left(
\begin{array}{ccc}
 \a & 0 & 0 \\
 0 & \b & 0 \\
 0 & 0 & 1 \\
\end{array}
\right)\mid \a, \b\in k^\ast\}
$$
acts as
$$
\left(
\begin{array}{ccc}
 \a & 0 & 0 \\
 0 & \b & 0 \\
 0 & 0 & 1 \\
\end{array}
\right)
.x_m=
\left(
\begin{array}{ccc}
 0 & 0 & \a\frac{a b c}{m} \\
 0 & -m & -\b\frac{(a+m) (b+m) (c+m)}{m} \\
 \a^{-1} & \b^{-1} & a+b+c+m \\
\end{array}
\right)
$$
Hence
$$
T.x_m\cap{\cal F}=\{x_m\}$$
The characteristic polynomial of $x_m(a,b,c)$ is $(X-a)(X-b)(X-c)$. Moreover, $\dim B.x_m(a,b,c)=5$, so that $\dim \O_{x_m(a,b,c)}=6$. We have shown that $x_m(a,b,c)$ is regular for every choice of $a$, $b$, $c \in k^\ast$. Now let $g$ be a regular element of $GL(3)$. Since $\O_g$ is determined by the characteristic polynomial of $g$, there exists $a$, $b$, $c\in k^\ast$ such that $x_m(a,b,c)\in \O_g$ for every $m\in k^\ast$. We take $x_m=x_m(a,b,c)$ for $m\in k^\ast$: the set ${\cal F}=\{x_m\mid m\in k^\ast\}$ is the required set.
\cvd

\begin{lemma} Let $H=Sp(4)$, any characteristic, $g$ a mixed regular element of $H$. Then there exists a subset
${\cal F}=\{x_m\mid m\in k^\ast\}$ of  $\O_g$ such that $(B(H).x_m)_{ m\in k^\ast}$ consists of pairwise distinct $B(H)$-orbits.
\end{lemma}

\pf
Without loss of generality we may assume 
$$
g=
\left(
\begin{array}{cccc}
 a & 0 & 0 & 0 \\
 0 & 1 & 0 & 1 \\
 0 & 0 & \frac{1}{a} & 0 \\
 0 & 0 & 0 & 1 \\
\end{array}
\right)
$$
with $a\not=\pm 1$. We put
$$
x_m=
\left(
\begin{array}{cccc}
 0 & 0 & -\frac{1}{m} & 0 \\
 0 & 0 & -1 & 1 \\
 m & m & \frac{a^2+m+1}{a} & \frac{m (-2 a+m+1)}{a} \\
 0 & -1 & -\frac{1}{a} & 2-\frac{m}{a} \\
\end{array}
\right)=
$$
$$
\left(
\begin{array}{cccc}
 0 & 0 & -1 & 0 \\
 0 & 0 & 0 & -1 \\
 1 & 0 & 0 & 0 \\
 0 & 1 & 0 & 0 \\
\end{array}
\right)
\left(
\begin{array}{cccc}
 m & 0 & 0 & 0 \\
 0 & -1 & 0 & 0 \\
 0 & 0 & \frac{1}{m} & 0 \\
 0 & 0 & 0 & -1 \\
\end{array}
\right)
\left(
\begin{array}{cccc}
 1 & 1 & \frac{a^2+m+1}{a m} & \frac{-2 a+m+1}{a} \\
 0 & 1 & \frac{1}{a} & \frac{m}{a}-2 \\
 0 & 0 & 1 & 0 \\
 0 & 0 & -1 & 1 \\
\end{array}
\right)\in w_0B
$$
for $m\in k^\ast$. Then $x_m$ is $H$-conjugate to $g$. In fact the characteristic polynomial of $x_m$ is \mbox{$(X-1)^2(X-a)(X-\frac1a)$}, and the $1$-eigenspace has dimension 1. 
Suppose $x_m$, $x_{m'}$ are $B$-conjugate. Then  $x_m$, $x_{m'}$ are $T$-conjugate and,
from a direct calculation, it follows that $T.x_m\cap{\cal F}=\{x_m\}$, hence $m=m'$.
\cvd

\end{document}